 \documentclass[final,3p,times]{elsarticle}

 \usepackage{graphics}

\usepackage{amssymb}
\usepackage{amsthm}

\usepackage{MnSymbol}

\usepackage{lineno}
\usepackage{color}                    
\newcommand\tstrut{\rule{0pt}{2.4ex}}

\biboptions{semicolon, sort&compress}
\usepackage{hyperref}
\hypersetup{colorlinks=true}
\usepackage{float}

\usepackage{color}
\newtheorem{theorem}{Theorem}
\newtheorem{lemma}{Lemma}
\newtheorem{corollary}{Corollary }
\newtheorem{remark}{Remark}
\newtheorem{definition}{Definition}

\newtheorem{question}{Question}

\newtheorem{fact}{Fact}

\theoremstyle{fact}
\newtheorem{examplex}[fact]{Example}
\newtheorem*{continuedex}{Example \continuedexref\space (Continued)}
\newtheorem*{example*}{Example}

\newcommand{\myqedsymbol}{\ensuremath{\triangle}}%

\newenvironment{examcont}[1]
{\newcommand{\continuedexref}{\ref*{#1}}\continuedex}
{\hfill\myqedsymbol\endcontinuedex}

\newenvironment{example}
{\pushQED{\qed}\examplex}
{\popQED\endexamplex}

\usepackage{graphicx}
\usepackage{caption}
\usepackage{subcaption}

\journal{Linear Algebra and its Applications}

\begin{document}

\begin{frontmatter}

\title{On Eigenvalues of Laplacian Matrix for a Class of Directed Signed Graphs}

\author[label1]{Saeed Ahmadizadeh}
\address[label1]{Department of Electrical and Electronic Engineering, Melbourne Information, Decision, and Autonomous Systems (MIDAS) Lab, The University of Melbourne, Parkville, VIC 3010, Australia}


\ead{sahmadizadeh@student.unimelb.edu.au}

\author[label1]{Iman Shames}
\ead{iman.shames@unimelb.edu.au}

\author[label2]{Samuel Martin}
\address[label2]{Universit\'{e} de Lorraine and CNRS, CRAN, Vandoeuvre-l\`{e}s-Nancy, France}
\ead{samuel.martin@univ-lorraine.fr}

\author[label1]{Dragan Ne\v{s}i\'{c}}
\ead{dnesic@unimelb.edu.au}

\begin{abstract}
The eigenvalues of the Laplacian matrix for a class of directed graphs with both positive and negative weights are studied.  First, a class of directed signed graphs is investigated in which one pair of nodes (either connected or not) is perturbed with negative weights. A necessary  condition is proposed to attain the following objective for the perturbed graph: the real parts of the non-zero eigenvalues of its Laplacian matrix are positive. A sufficient condition is also presented that ensures the aforementioned objective for unperturbed graph. It is then highlighted the case where the condition becomes necessary and sufficient. Secondly, for directed graphs, a subset of pairs of nodes are identified where if any of the pairs is connected by an edge with infinitesimal negative weight, the resulting Laplacian matrix will have at least one eigenvalue with negative real part. Illustrative examples are presented to show the applicability of our results.

\end{abstract}


\end{frontmatter}



\section{Introduction}
\label{sec:introduction}

Complex networks arise in a wide range of applications such as social networks and multi-agent systems.  The analysis of such networks is often split into two parts. In the first part, the dynamical properties of each agent such as passivity or dissipativity~\cite{arcak2007passivity,stan2007analysis} are investigated. These properties generally facilitate the analysis in the second part in which every agent is treated as a node in a graph, and their interconnections are modeled as the edges.

The evolution of the states of each agent is primarily influenced by the relative information from its neighbors This type of interconnection is mathematically described by the Laplacian matrix which captures the structure of the network graph. In general, this graph can be directed with both positive and negative weights. As an example, in the network of neurons, the connections between presynaptic neurons and postsynaptic ones are directed. Furthermore, the coupling weights between excitatory neurons to other excitatory neurons are positive while the coupling weights between inhibitory neurons to other inhibitory neurons are negative~\cite{hoppensteadt2012weakly,rodriguez1999perception}.

There is a wide range of applications where the eigenvalues of the Laplacian matrix play a crucial role in the behavior of the network. For instance, it is well-known that consensus can be reached if the Laplacian matrix is allowed to have a single zero eigenvalue with all non-zero eigenvalues having positive real parts~\cite{olfati2007consensus}. Synchronization in the network of linear and nonlinear systems highly depends on the structure of network which itself determines the eigenvalues of the Laplacian matrix~\cite{belykh2006generalized}. In most studies, synchronization cannot be guaranteed if the Laplacian matrix has an eigenvalue with negative real part~\cite{scardovi2009synchronization,shafi2013synchronization,ahmadizadeh2016synchronization}. The existence of both positive and negative weights also leads to clustering in the network that can be demonstrated in terms of eigenvectors and eigenvalues of the Laplacian matrix~\cite{shames2014manipulating,xia2011clustering}.

The spectral characterization of Laplacian matrix has been a subject of many research activities that is well-understood for  undirected graphs with non-negative weights~\cite{chung1997spectral}. For directed graphs with positive weights, there is a relation between the so-called normalized Laplacian matrix and the stochastic matrices that is employed to explore the locations of eigenvalues of the Laplacian matrix~\cite{agaev2005spectra,caughman2006kernels,brualdi2010spectra}. In the presence of negative weights, this relation is hard to establish and, consequently, the spectral characterization of Laplacian matrix becomes  more challenging. For undirected graphs, there exist powerful results that mainly provide bounds on the number of negative and positive eigenvalues. ~\cite{bauer2012normalized,bronski2014spectral,bronski2015spectral}. Recently, a connection between the number of positive eigenvalues of the Laplacian matrix and negative weights has been developed for undirected graphs in~\cite{bronski2016graph}. More specifically, it has been proven that the number of positive eigenvalues equals the number of negative weights in the graph minus the number of positive eigenvalues of the associated  cycle intersection matrix~\cite[Theorem 2.9]{bronski2016graph}. However, no conditions on the magnitude of the negative weights has been identified to ensure the existence of eigenvalues with negative real parts.

In many applications, e.g.~see~\cite{trentelman2013robust}, the objective is to prevent the Laplacian matrix from having eigenvalues with negative real parts. Robustness of uncertain undirected networks has been recently studied where the negative weight has been incorporated as uncertainty in the network. Under certain assumptions on the distribution of negative weights, the robustness of network has also been analyzed in the presence of multiple negative weights~\cite{zelazo2015robustness}. These results have recently been extended to a general case~(undirected graphs with arbitrary numbers of negative weights) using tools from electrical circuit theory~\cite{chen2016characterizing}. It has been argued that these results can be interpreted using the notion of effective resistance originally introduced in electrical networks~\cite{klein1993resistance}. The robustness of the consensus algorithm over directed signed graphs has been studied in~\cite{mukherjee2016consensus}.

Our results concern analyzing the eigenvalues of the Laplacian matrix for directed signed graphs. Similar to the undirected graphs, we obtain a necessary/sufficient conditions~(for some cases necessary and sufficient) which provides an upper bound of negative weights between any pairs of nodes which guarantees that none of the eigenvalues of the Laplacian matrix has a negative real part. It was not possible to employ the machinery employed in~\cite{chen2016characterizing,zelazo2015robustness} for undirected graphs to derive the results of this paper as the Laplacian matrix of the directed graphs are not generally symmetric. The employed methodology in~\cite{zelazo2015robustness} has been recently extended to deal with directed graphs~\cite{mukherjee2016consensus} where sufficient conditions for the  upper bound on a single negative weight has been derived via Nyquist stability criteria.  The necessary result of our paper considers the case in which both edges between any arbitrary pairs of nodes are perturbed with  negative weights.  Our sufficiency result is more general than the main result of~\cite{mukherjee2016consensus} since we also allow  perturbing two edges between two nodes with the same negative weight to the signed directed graph.  Our results cover a more general set of graphs as a graph with multiple negative edges might satisfy the assumption of the theorem, while~\cite[Theorem 1]{mukherjee2016consensus} only is applied to graphs with no negative edges. Even though the results of our paper are interpreted via Nyquist criteria, our approaches are different from~\cite{mukherjee2016consensus}. We also highlight the case where the condition becomes necessary and sufficient.  Furthermore,  it is argued that for directed graphs, the recently proposed notion of effective resistance~\cite{young2015new} is not applicable to interpret the obtained upper bound. By partitioning the nodes of the graph into some sets, we identify ``sensitive pairs of nodes" with the following property: If there exists at least one edge with sufficiently small negative weight, the Laplacian matrix has at least one eigenvalue with negative real part. This result is different from the main result in~\cite[Theorem 2.10]{bronski2014spectral} which established a lower bound and upper bound for the number of negative eigenvalues.

This paper is organized as follows. In Section~\ref{sec:preliminaries}, we introduce the required notation and preliminaries. Section~\ref{sec:problem_results} states the underlying questions studied in the current manuscript following by the main results, illustrative examples, and discussions. Section~\ref{sec:application} demonstrates the applicability of this work through consensus in social networks. Conclusion is presented in the last section and the proofs of some lemmas are included in Appendix.

\section{Preliminaries and Graph Notation}
\label{sec:preliminaries}

Throughout this paper, $I_N \in \mathbb{R}^
{N \times N}$, $\mathbf{1}_N \in \mathbb{R}^
{N}$, and $\mathbf{0}_N \in \mathbb{R}^
{N}$ denote the ${N \times N}$ identity matrix, the $N$-dimensional vectors containing $1$, and $0$ in every entry, respectively. The standard bases in $\mathbb{R}^N$ are represented by $\{ \mathbf{e}_1,\dots,\mathbf{e}_N \}$ where $\mathbf{e}_i$ is the $i^{th}$  column of  $I_N$. The $2$-norm of a vector $x \in \mathbb{R}^N$ is shown by $ \| x \|$. For a complex
variable, vector or matrix, $\Re(\cdot)$ and $\Im(\cdot)$ stand for the real and imaginary parts. For a matrix $A \in \mathbb{R}^{N \times N}$, $\text{Spec}(A) = \{ \lambda_i(A) \}_{i=1}^{N}$ denotes the set of eigenvalues of $A$ where  $\Re (\lambda_1) \leq \Re (\lambda_2) \leq  \dots \leq \Re(\lambda_N)$. An eigenvalue $\lambda_i(A)$ is called semisimple if its algebraic and geometric multiplicities are equal.
The operator $\text{diag}(\cdot)$ constructs a block diagonal matrix from its arguments. The entry in the $i^{th}$ row and $j^{th}$ column of a matrix $A$ is represented by $ [A ]_{ij}$, while the $i^{th}$ entry of a vector $x$ is denoted by $[x]_i$. For a set $\mathcal{A}$, its cardinality is denoted by $|\mathcal{A}|$.

A weighted directed signed graph $\mathcal{G}$ is represented by the triple $\mathcal{G}(\mathcal{V},\mathcal{E}, \mathcal{W})$ where $\mathcal{V}=\{ 1,\dots,N \}$ is the nodes set, $\mathcal{E} \subset \mathcal{V} \times \mathcal{V}$ is the edge set, and $\mathcal{W} :  \mathcal{V} \times \mathcal{V} \to \mathbb{R}$ is a weight function that maps each $(i,j) \in \mathcal{E}$ to a nonzero scalar $a_{ij}$ and returns 0 for all other $(i,j)\not\in\mathcal{E}$. The adjacency matrix $ A \in \mathbb{R}^{N \times N}$ captures the interconnection between the nodes in the graph where $[A]_{ij} = a_{ij} \neq 0$ iff $(i,j) \in \mathcal{E}$. For the edge $(i,j)$, we follow the definition corresponding to a sensing convention which indicates that node $i$ receives information form node $j$ or equivalently, the node $j$ influences the node $i$; see~\cite{young2015new} for more information. For each node $i \in \mathcal{V}$, $\mathcal{N}_i$ denotes the set of its neighbors, i.e., $\mathcal{N}(i)= \{j~|~ a_{ij} \neq 0 \}$.

For a given graph $\mathcal{G}(\mathcal{V},\mathcal{E},\mathcal{W})$ and a set $\overline{\mathcal{V}} \subseteq \mathcal{V}$, the corresponding induced subgraph is denoted by $\mathcal{G}(\overline{\mathcal{V}},\overline{\mathcal{E}},\overline{\mathcal{W}})$, where the set $\overline{\mathcal{E}}$ is defined as $\overline{\mathcal{E}} = \{ (i,j) \in \mathcal{E}~| i,j \in \overline{\mathcal{V}} \}$, and $\overline{\mathcal{W}} : \overline{\mathcal{V}} \times \overline{\mathcal{V}} \to \mathbb{R}$ is defined as $\overline{\mathcal{W}} (i,j) = \mathcal{W}(i,j)$.
In order to categorize edges in terms of the sign of their values,  we define the sets  $\mathcal{E}^+ = \{ (i,j)\,|\,a_{ij} > 0 \}$, and $\mathcal{E}^{-}=\mathcal{E} \backslash \mathcal{E}^+ =\{ (i,j)\,|\,a_{ij} <0 \}$. We call the edges in $\mathcal{E}^+$ and $\mathcal{E}^-$ positive edges and negative edges, respectively. Subsequently, for a signed graph $\mathcal{G}(\mathcal{V},\mathcal{E}, \mathcal{W})$, we denote the subgraph with non-negative weights by $\mathcal{G} (\mathcal{V},\mathcal{E}^+,\mathcal{W}^+)$ where $\mathcal{W}^+ : \mathcal{V} \times \mathcal{V} \to \mathbb{R}_{ \geq  0}$ is defined as $\mathcal{W}^+ (i,j) = \mathcal{W}(i,j)$ for all $ (i,j)\in \mathcal{E}^+ $ and $\mathcal{W}^+ (i,j) =0$ for all $ (i,j)\nin \mathcal{E}^+ $. Similarly, for a signed graph $\mathcal{G}(\mathcal{V},\mathcal{E}, \mathcal{W})$, we denote the subgraph with non-positive weights by $\mathcal{G} (\mathcal{V},\mathcal{E}^-,\mathcal{W}^-)$.
The superposition of two signed directed graphs $ \mathcal{G}_1(\mathcal{V},\mathcal{E}_1, \mathcal{W}_1) \bigoplus \mathcal{G}_2(\mathcal{V},\mathcal{E}_2, \mathcal{W}_2)$ is a new graph $\mathcal{G}(\mathcal{V},\mathcal{E}, \mathcal{W})$ where $\mathcal{E}= \mathcal{E}_1 \cup \mathcal{E}_2$ and, $\mathcal{W}(i,j)= \mathcal{W}_1(i,j) + \mathcal{W}_2(i,j)$ for every $(i,j) \in \{ \mathcal{V} \times \mathcal{V} \}$.

\renewcommand{\thefootnote}{\fnsymbol{footnote}}

For a directed graph, the in-degree and out-degree of node $i$ are defined as $d^{in}_i=\sum\limits_{j} a_{ji} $ and $d^{out}_i=\sum\limits_{j} a_{ij}$ respectively.
The Laplacian matrix\footnote{The current definition of Laplacian matrix has been inspired from a large variety of applications in control such as consensus~\cite{zelazo2015robustness}, security analysis of complex networks~\cite{teixeira2010networked,shames2011distributed}, and synchronization in networks of oscillators~\cite{shafi2013synchronization, ahmadizadeh2016synchronization}. However, there is another way to define the Laplacian matrix for weighted signed graphs in which the in-degree and out-degree of node $i$ are defined as $d^{in}_i=\sum\limits_{j} | a_{ji} | $ and $d^{out}_i=\sum\limits_{j} | a_{ij} |$ respectively~\cite{dong2016laplacian}.    } $L \in \mathbb{R}^{N \times N}$ is defined by $L=D-A$ where $D= \text{diag} \{ d^{out}_1,\dots,d^{out}_N \}$. Since the rows of the Laplacian matrix add to zero, $\mathbf{1}_N$ is always one of its eigenvector that corresponds to the eigenvalue $0$. This eigenvalue is called the trivial eigenvalue while the rest of eigenvalues is called non-trivial eigenvalues.

Let $\Pi=I_N - \frac{1}{N} \mathbf{1}_N \mathbf{1}_N^T$ denote the orthogonal projection matrix onto the subspace of $\mathbb{R}^N$ perpendicular to $\mathbf{1}_N$. The matrix $\Pi$ is symmetric and since $L \mathbf{1}_N = \mathbf{0}$, $L \Pi = L$ and $ \Pi L^T = L^T$ for any graph.
We define a matrix $Q \in \mathbb{R}^{(N-1) \times N}$ whose rows are the orthonormal bases for $\text{span} \{\mathbf{1}_N \}^{\perp}$ where $\perp$ denotes the orthogonal complement of the space. Hence, $Q^T$ is a full column rank matrix. On $\text{span} \{\mathbf{1}_N \}^{\perp}$, the Laplacian matrix is equivalent to the so-called \textit{reduced Laplacian} $\bar{L} \in \mathbb{R}^{(N-1) \times (N-1)}$ which is defined by~\cite{young2015new},
\begin{equation}
	\label{eq:reduced_laplacian}
	\begin{aligned}
		\bar{L} := Q L Q^T.
	\end{aligned}
\end{equation}

A \emph{path} of length $r$ from $i_1\in\mathcal{V}$ to $i_r\in\mathcal{V}$ in graph $\mathcal{G}$ is a sequence $(i_1,i_2,\dots,i_r)$ of distinct nodes in $\mathcal{V}$ where $i_{j+1}$ is a neighbor of $i_j$ for all $j = 1,\dots,r-1$. If there exists a path (no path) from the node $j$ to the node $i$, then the node $i$ is (not) reachable from node $j$. We use $j \rightarrowtail i $ ($j \nrightarrowtail i $) to show the existence (absence) of path from $ j$ to $i$.  A node $i$ is a globally reachable node if it is reachable from all other nodes of the graph.
Similar to~\cite{young2015new}, we say that two nodes $i$ and $j$ are connected if the graph contains two paths, one starting from node $i$ and the other one from $j$ that both terminates at the same node. A graph $\mathcal{G}$ is connected if every pair of nodes is connected. This notion of a connected graph corresponds to the scrambling matrices~\cite{seneta1979coefficients}. It has been shown that the graph is connected if and only if there exists at least one \textit{globally reachable node}; the node, to which, there exists at least one path from every node in the graph~\cite{young2015new}. A graph $\mathcal{G}$ is strongly connected if for every $i \in \mathcal{V}$ and $j \in \mathcal{V}$, $ i \rightarrowtail j$. Hence, the graph $\mathcal{G}$ is strongly connected if and only if every node of the graph is a globally reachable node.

For every node $i \in \mathcal{V}$, the reachable set $\mathcal{R}(i)$ is the union of the node $i$ and all nodes from which there exists a path to node $i$, i.e. $\mathcal{R}(i) = \{ i \} \cup \{ j \in \mathcal{V}~|~j \rightarrowtail i  \}$. Next the reachable, exclusive and common sets associated with a graph $\mathcal{G}$ are defined.

\begin{definition}
	\label{def:reachable_sets}
	For a given graph $\mathcal{G}(\mathcal{V},\mathcal{E},\mathcal{W})$, a set of nodes $\mathcal{R} \subseteq \mathcal{V}$ is called a reach set if (a) there exists at least one node $ i \in \mathcal{V}$ such that $\mathcal{R}(i) = \mathcal{R}$, and (b) it is maximal, i.e. there is no $j \in \mathcal{V}$ such that $\mathcal{R}(i) \subset \mathcal{R}(j)$~(properly). Let  $\mathcal{R}_1,\dots,\mathcal{R}_d$ denote all reaches of the graph $\mathcal{G}$. We associate a set of reaching nodes $\mathcal{U}_k$ with each reach set $\mathcal{R}_k$ and it is defined as $\mathcal{U}_k = \{ i \in \mathcal{V}~|~\mathcal{R}(i)=\mathcal{R}_k \}$.
	For each set $\mathcal{R}_k$, the exclusive set $\mathcal{X}_k$ and the common set $\mathcal{C}_k$ are defined as $\mathcal{X}_k= \mathcal{R}_k \backslash \cup_{i \neq k} \mathcal{R}_i  $ and $\mathcal{C}_k= \mathcal{R}_k \backslash \mathcal{X}_k$, respectively.
\end{definition}

\begin{remark}
	\label{remark:different_convention}
	The description of the sets in Definition~\ref{def:reachable_sets} differs from~\cite[Definition 2.6]{caughman2006kernels} as we use the sensing convention rather the information flow one. For a given graph with sensing convention, the same results can be obtained by using the information flow convention if the direction of each edge is reversed.
\end{remark}

To obtain an intuition behind the concepts in Definition~\ref{def:reachable_sets}, we see the weight $a_{ij} \neq 0$ as an influence weight from $j$ to $i$. In this case, a reach set $\mathcal{R}_k$ is a set in which at least one of node influences all others directly or indirectly, and this set cannot be included in a bigger such set. Set $\mathcal{U}_k$ is the set of nodes in $\mathcal{R}_k$ which influences all others in $\mathcal{R}_k$. Set $\mathcal{X}_k$ is the set of nodes in $\mathcal{R}_k$ which are in no other reach set $\mathcal{R}_q$ with $q \neq k$. Set $\mathcal{C}_k$ is the set of nodes in $\mathcal{R}_k$ which are also in some other reach set $\mathcal{R}_q$ with $q \neq k$. The following example shows how these sets are defined for a given graph.

\begin{example}
	\label{exa:EXS1}
	The graph shown in Figure~\ref{fig:fig1} has three reach sets  $\mathcal{R}_1 = \{ 1,2,8,9,10,11,12\}$,~$\mathcal{R}_2 = \{ 3,4,5,8,9,10,11,12\}$ and $\mathcal{R}_3 = \{ 6,7\}$. By using the definitions of $\mathcal{X}_k$ and $\mathcal{C}_k$ for $k=1,2,3$, the exclusive sets and common sets includes $\mathcal{X}_1 = \{ 1,2\}$, $\mathcal{X}_2 = \{ 3 ,4 ,5 \}$, $\mathcal{X}_3 = \{ 6, 7 \}$, $\mathcal{C}_1 = \mathcal{C}_2 = \{ 8,9,10,11,12 \}$ and $\mathcal{C}_3 = \emptyset$. The reaching sets are described by $\mathcal{U}_1 = \{ 1,2 \}$, $\mathcal{U}_2 = \{ 3 \}$ and $\mathcal{U}_3 = \{ 7 \}$.
	\begin{figure}[ht]
		\centering
		\includegraphics[scale=0.6]{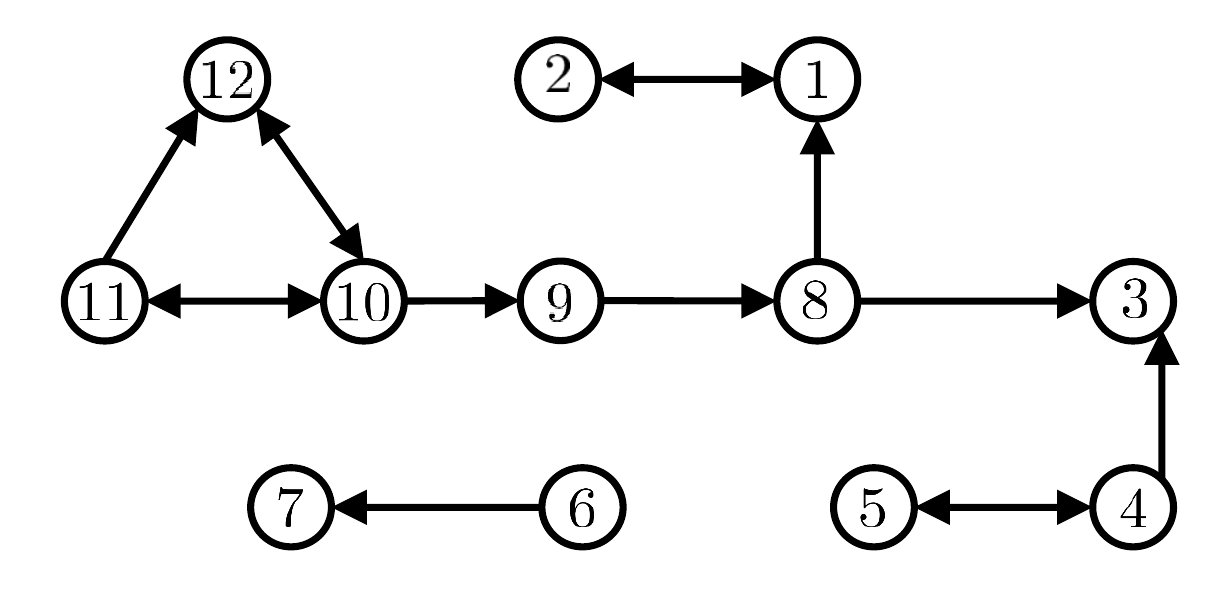}
		\caption{ A directed graph with $12$ nodes in Example~1.}
		\label{fig:fig1}
	\end{figure}
\end{example}

In the following lemma, some properties of these sets are presented.

\begin{lemma}
	\label{lemma:sets_properties}
	Given a graph $\mathcal{G}(\mathcal{V},\mathcal{E},\mathcal{W})$ and the corresponding sets in Definition~\ref{def:reachable_sets}, the following properties hold for every $k=1,\dots,d$.
	\begin{enumerate}
		\item { $\mathcal{X}_k \cap \mathcal{C}_k = \emptyset$.
		}
		\item { $\mathcal{U}_k \subseteq \mathcal{X}_k$.
		}
		\item{ For every $i \in \mathcal{U}_k$, $j \in \mathcal{X}_k \backslash \mathcal{U}_k$, $p \in \mathcal{X}_k$, and $m \in \mathcal{C}_k$, we have
			\begin{enumerate}
				\item {$(i,j) \notin \mathcal{E}$.
				}
				\item {$(p,m) \notin \mathcal{E}$.
				}
			\end{enumerate}
		}
		\item{ For every $\mathcal{U}_k$, the corresponding induced subgraph $\mathcal{G}(\mathcal{U}_k,\mathcal{E}_{\mathcal{U}_k},\mathcal{W}_{\mathcal{U}_k})$ is strongly connected.
		}
		\item{If $\mathcal{C}_k = \emptyset$ for some $k$, then $\mathcal{R}_k = \mathcal{X}_k$ which means that induced subgraph~$\mathcal{G}(\mathcal{R}_k,\mathcal{E}_{\mathcal{R}_k},\mathcal{W}_{\mathcal{R}_k})$ is disconnected from the rest of the graph.
		}
	\end{enumerate}
\end{lemma}
\textit{Proof of Lemma~\ref{lemma:sets_properties}:}
\textit{Proof of 1.}
This property is a direct consequence of how the sets $\mathcal{C}_k$, $\mathcal{X}_k$ are defined.

\textit{Proof of 2.} Choose an arbitrary node $i \in \mathcal{U}_k$. From the definition of $\mathcal{U}_k$ it follows that $\mathcal{R}(i) = \mathcal{R}_k$ and by definition of $\mathcal{R}(i)$ we have that $i \in \mathcal{R}_k$. We complete the proof by contradiction. Assume $i \in \mathcal{C}_k$ which means there exists at least one $q \neq k$ such that $ i \in \mathcal{C}_q$ and subsequently $ i \in \mathcal{R}_q$. This means $\mathcal{R}_k \subset \mathcal{R}_q$ which contradicts maximality of~$\mathcal{R}_k$. Hence $ i \notin \mathcal{C}_k$.

\textit{Proof of 3(a).} Since $i \in \mathcal{U}_k$, then $\forall u \in \mathcal{R}_k, u  \rightarrowtail i$. For the purpose of showing contradiction, suppose that $(i,j) \in \mathcal{E}$. Then $ i \rightarrowtail j$ and subsequently $\forall u \in \mathcal{R}_k, u  \rightarrowtail j$ which is equivalent to  $j \in \mathcal{U}_k$. This contradicts $j \in \mathcal{X}_k \backslash \mathcal{U}_k$.

\textit{Proof of 3(b).} Assume that $p \in \mathcal{X}_k$ and $m \in \mathcal{C}_k$. For the purpose of showing contradiction, assume $(p,m) \in \mathcal{E}$ which is equivalent to $ j \rightarrowtail m$. Since $m \in \mathcal{C}_k$, there exists at least $q$ such that $m \in  \mathcal{R}_q$ which along with the existence of $ j \rightarrowtail m$ implies $j \in \mathcal{R}_q$. This concludes $\mathcal{X}_k \cap \mathcal{R}_q \neq \emptyset$ which is in contradiction with the definition of $\mathcal{X}_k$.

\textit{Proof of 4.} We prove the fourth property by contradiction. Since a graph with single node is strongly connected, without loss of generality, we assume $|\mathcal{U}_k| > 1$. Assume that $\mathcal{G}_k$ is not strongly connected. Hence, there exists at least two nodes $i,~j \in \mathcal{U}_k$ such that $i \nrightarrowtail j$. This means $i \notin \mathcal{R}(j) $ and since by definition of $\mathcal{R}(i),~i\in \mathcal{R}(i)$, we have  $\mathcal{R}(i) \neq \mathcal{R}(j)$. This contradicts the assumption $\mathcal{R}(i) = \mathcal{R}(j) = \mathcal{R}_k$. Hence, $\mathcal{G}_k$ is strongly connected.

\textit{Proof of 5.} If $\mathcal{C}_k = \emptyset$, from Definition~\ref{def:reachable_sets}, we then obtain  $\mathcal{R}_k \backslash \mathcal{X}_k = \emptyset \Longleftrightarrow \mathcal{R}_k = \mathcal{X}_k$. Furthermore,  $\mathcal{R}_k \cap \mathcal{R}_i = \emptyset$ for all $i \neq k$ which means $$\forall u \in \mathcal{R}_k , \forall v \in \mathcal{R}_i~\Longrightarrow (u,v) \notin \mathcal{E}~\text{and}~(v,u) \notin \mathcal{E}. $$
This completes the proof of the fifth property. \qed

We now turn to Example 1 to illustrate that Lemma\ref{lemma:sets_properties} holds in this example.

\begin{examcont}{exa:EXS1}
	It is straightforward to check that the properties $1$ and $2$ in Lemma~\ref{lemma:sets_properties} hold. To check the third property, it is observed that $(3,4) \notin \mathcal{E}$ as $3 \in \mathcal{U}_2$ and $4 \in \mathcal{X}_2 \backslash \mathcal{U}_2$. Note that the node $4$ cannot be in the set $\mathcal{U}_2$ as $\mathcal{R}(4) = \{ 3, 4, 5\} \subset \mathcal{R}_2$.  It is observed that all subgraphs induced by $\mathcal{U}_k,~k=1,2,3$ are strongly connected. We note that $\mathcal{C}_3$ is empty which means the original graph has a disconnected subgraph.
\end{examcont}

For a graph with non-negative weights, we list the properties of its Laplacian matrix and its eigenvectors in the following lemma.

\begin{lemma}
	\label{lemma:eigenvectors_properties}
	Consider a graph $\mathcal{G}(\mathcal{V},\mathcal{E},\mathcal{W})$ with non-negative weights and Laplacian matrix $L$. Assume that the graph has $d$\footnote{The number of reach sets $d$ can be found by identifying the reachable set of each node and then checking conditions (a) and (b) in Definition~\ref{def:reachable_sets}.}
	reach sets~$\mathcal{R}_k$ with the corresponding sets $\mathcal{X}_k$, $\mathcal{U}_k$, and $\mathcal{R}_k$ according to Definition~\ref{def:reachable_sets}. Then, zero is a semisimple eigenvalue of  $L$ with multiplicity $d$.
	Denote the right and the left eigenvectors associated with the zero eigenvalue of $L$ that form the orthogonal bases for its eigenspace by $\gamma_1,\dots,\gamma_d \in \mathbb{R}^N$ and $\mu_1,\dots,\mu_d \in \mathbb{R}^N$, respectively. Each $\gamma_k,~k=1,\dots,d$ can be chosen to satisfy the following conditions:
	\begin{enumerate}
		\item {
			$[\gamma_k]_i = 1$ for $i \in \mathcal{X}_k$;
		}
		\item {$[\gamma_k]_i = 0$ for $i \in \mathcal{V} \backslash \mathcal{R}_k$;
		}
		\item{$[\gamma_k]_i \in (0,1) $ for $i \in   \mathcal{C}_k$;
		}
		\item{$\sum\limits_{k=1}^{d} \gamma_k = \mathbf{1}_N$.
		}
	\end{enumerate}
	Furthermore, each $\mu_k$ can be chosen such that their entries satisfy the following conditions:
	\begin{enumerate}
		\item {
			$[\mu_k]_i \in (0,1)  $ for $i \in \mathcal{U}_k$;
		}
		\item {$[\mu_k]_i = 0$ for $i \in \mathcal{V} \backslash \mathcal{U}_k$;
		}
		\item{$\sum\limits_{i \in \mathcal{U}_k}^{}[\mu_k]_i=1$
		}
	\end{enumerate}
\end{lemma}
\textit{Proof.}
By taking into account Remark~\ref{remark:different_convention}, the fact that the eigenvalue is semisimple and the results regarding $\gamma_k$ correspond to~\cite[Corollary~4.1]{caughman2006kernels}. Denote $A$ the adjacency matrix of $\mathcal{G}$. By relabeling the nodes, they can be partitioned such that the first  $|\mathcal{X}_1|$ nodes in the first partition belong to the set $\mathcal{X}_1$, the second $|\mathcal{X}_2|$ nodes in the second partition belong to the set $\mathcal{X}_2$, and so on, and the remaining nodes belong to the set $\mathcal{C} = \bigcup\limits_{k=1}^{d} \mathcal{C}_k$. This means that there exist maps~$\Pi_{\mathcal{V}} : \mathcal{V} \to \widetilde{\mathcal{V}}$, $\Pi_{\mathcal{E}} : \mathcal{E} \to \widetilde{\mathcal{E}}$ such that the new graph~$\widetilde{\mathcal{G}}(\widetilde{\mathcal{V}},\widetilde{\mathcal{E}},\widetilde{\mathcal{W}})$ has $d$ reach sets~$\widetilde{\mathcal{R}}_k = \Pi_{\mathcal{V}} \circ \mathcal{R}_k $, exclusive sets~$\widetilde{\mathcal{X}}_k = \Pi_{\mathcal{V}} \circ \mathcal{X}_k $, and reaching nodes sets~$\widetilde{\mathcal{U}}_k = \Pi_{\mathcal{V}} \circ \mathcal{U}_k $, where $\circ$ denotes the element-wise operation on the set. There is also a permutation matrix $P \in \mathbb{R}^{N \times N}$ that relates the adjacency and Laplacian matrices of $\widetilde{\mathcal{G}}$ to the ones of $\mathcal{G}$ via $\tilde{A} = PAP^{-1}$ and $\tilde{L} = PLP^{-1}$.

By employing Properties~$1-3$ in Lemma~\ref{lemma:sets_properties}, $\tilde{A}$ can be written as
\begin{equation}
	\label{eq:adjacency_general}
	\left[\begin{array}{cccc|c}
		\tilde{A}_{\widetilde{\mathcal{X}}_1 \widetilde{\mathcal{X}}_1}&\mathbf{0}&\hdots&\hdots&\mathbf{0}\\
		\mathbf{0}&\tilde{A}_{\widetilde{\mathcal{X}}_2 \widetilde{\mathcal{X}}_2}&\hdots&\hdots&\mathbf{0}\\
		\vdots&\ddots &\ddots&\hdots&\vdots\\
		\mathbf{0}&\mathbf{0}&\hdots&\tilde{A}_{\widetilde{\mathcal{X}}_d \widetilde{\mathcal{X}}_d}&\mathbf{0} \\   \hline
		\tilde{A}_{\widetilde{\mathcal{C}} \widetilde{\mathcal{X}}_1}&\tilde{A}_{\widetilde{\mathcal{C}} \widetilde{\mathcal{X}}_2}&\hdots&\tilde{A}_{\widetilde{\mathcal{C}} \widetilde{\mathcal{X}}_d}&\tilde{A}_{\widetilde{\mathcal{C}} \widetilde{\mathcal{C}}} \tstrut\\
	\end{array}\right],
\end{equation}
where $\tilde{A}_{S S'}$ means the matrix of weights induced by set $S$ and $S'$. Due to the property 3(a) in Lemma~\ref{lemma:sets_properties}, each block $\tilde{A}_{\widetilde{\mathcal{X}}_k \widetilde{\mathcal{X}}_k}$ is decomposed to
\begin{equation}
	\label{eq:adjacency_diagonal_block}
	\left[\begin{array}{cc}
		\tilde{A}_{\widetilde{\mathcal{U}}_i \widetilde{\mathcal{U}}_i}&\mathbf{0}\\
		\tilde{A}_{\widetilde{\mathcal{M}}_i \widetilde{\mathcal{U}}_i}&\tilde{A}_{\widetilde{\mathcal{M}}_i \widetilde{\mathcal{M}}_i}\\
	\end{array}\right],
\end{equation}
with $\widetilde{\mathcal{M}}_i = \widetilde{\mathcal{X}}_i \backslash \widetilde{\mathcal{U}}_i$. Consequently, the Laplacian matrix $\tilde{L}$ is represented by
\begin{equation}
	\label{eq:laplacian_general}
	\left[\begin{array}{cccc|c}
		\tilde{L}_{\widetilde{\mathcal{X}}_1 \widetilde{\mathcal{X}}_1}&\mathbf{0}&\hdots&\hdots&\mathbf{0}\\
		\mathbf{0}&\tilde{L}_{\widetilde{\mathcal{X}}_2 \widetilde{\mathcal{X}}_2}&\hdots&\hdots&\mathbf{0}\\
		\vdots&\ddots &\ddots&\hdots&\vdots\\
		\mathbf{0}&\mathbf{0}&\hdots&\tilde{L}_{\widetilde{\mathcal{X}}_d \widetilde{\mathcal{X}}_d}&\mathbf{0} \\   \hline
		\tilde{L}_{\widetilde{\mathcal{C}} \widetilde{\mathcal{X}}_1}&\tilde{L}_{\widetilde{\mathcal{C}} \widetilde{\mathcal{X}}_2}&\hdots&\tilde{L}_{\widetilde{\mathcal{C}} \widetilde{\mathcal{X}}_d}&\tilde{L}_{\widetilde{\mathcal{C}} \tilde{\mathcal{C}}} \tstrut\\
	\end{array}\right],
\end{equation}
where the block diagonal terms are
\begin{equation}
	\label{eq:Laplacian_diagonal_block}
	\left[\begin{array}{cc}
		\tilde{L}_{\widetilde{\mathcal{U}}_k \widetilde{\mathcal{U}}_k}&\mathbf{0}\\
		\tilde{L}_{\widetilde{\mathcal{M}}_k \widetilde{\mathcal{U}}_k}&\tilde{L}_{\widetilde{\mathcal{M}}_k \widetilde{\mathcal{M}}_k}\\
	\end{array}\right],
\end{equation}
and $\tilde{L}_{\widetilde{\mathcal{U}}_k \widetilde{\mathcal{U}}_k} \in \mathbb{R}^{|\widetilde{\mathcal{U}}_k| \times |\widetilde{\mathcal{U}}_k|}$ is the Laplacian matrix of the induced subgraph $\widetilde{\mathcal{G}}_k(\widetilde{\mathcal{U}}_k,\widetilde{\mathcal{E}}_{\widetilde{\mathcal{U}}_k},\widetilde{\mathcal{W}}_{\widetilde{\mathcal{U}}_k})$ which is strongly connected (property~4 in Lemma~\ref{lemma:sets_properties}).
Hence, there exists a left eigenvector $\tilde{\nu}_k \in \mathbb{R}^{|\widetilde{\mathcal{U}}_k|} $ that corresponds to the unique zero eigenvalue of $\tilde{L}_{\widetilde{\mathcal{U}}_k \widetilde{\mathcal{U}}_k}$, i.e. $\tilde{\nu}_k^T \tilde{L}_{\widetilde{\mathcal{U}}_k \widetilde{\mathcal{U}}_k} = \mathbf{0}^T$, with the following properties~\cite{brualdi2010spectra},
\begin{enumerate}
	\item{
		$[\tilde{\nu}_k]_j \in (0,1)$ for $i=1,\dots, |\widetilde{\mathcal{U}}_k|$.
	}
	\item{$\sum\limits_{i =1}^{|\widetilde{\mathcal{U}}_k|}[\tilde{\nu}_k]_i=1$
	}
\end{enumerate}
In light of the discussion above, it is observed that
\begin{equation}
	\label{eq:mu_k}
	\tilde{\mu}_k = [\mathbf{0}^T_{|\widetilde{\mathcal{X}}_1|}~\mathbf{0}^T_{|\widetilde{\mathcal{X}}_2|}~\dots~\tilde{\nu}_k^T \mathbf{0}^T_{|\widetilde{\mathcal{M}}_k|}~\dots~\mathbf{0}^T_{|\widetilde{\mathcal{X}}_d|}~\mathbf{0}^T_{|\widetilde{\mathcal{C}}|}]^T,
\end{equation}
is a left eigenvector for $\tilde{L}$ corresponding to the zero eigenvalue. Since $\widetilde{\mathcal{U}}_i \cap \widetilde{\mathcal{U}}_j = \emptyset$ for $i \neq j$, $\tilde{\mu}_j^T \tilde{\mu}_i = 0$ which means $\tilde{\mu}_1,\dots,\tilde{\mu}_d$ are orthogonal bases for the eigenspace of zero eigenvalues. As $\tilde{L} = P L P^{-1}$, we have $\mu_k = P^T \tilde{\mu}_k$ which along with ~(\ref{eq:mu_k}) and the properties of $\tilde{\nu}_k$ lead to Properties $1-3$ of $\mu_k$. This completes the proof. \qed

\begin{examcont}{exa:EXS1}
	Now, we illustrate the results of Lemma~\ref{lemma:eigenvectors_properties} by assigning the following weights to the edges of the graph in Figure~\ref{fig:fig1}: $a_{12}=2$, $a_{21}=1$, $a_{43}=3$, $a_{45}=5$, $a_{54}=4$, $a_{67}=7$, $a_{81}=1$, $a_{83}=3$,  $a_{98}=8$, $a_{10,9}=9$, $a_{10,11}=11$, $a_{10,12}=12$, $a_{11,10}=10$, $a_{11,12}=12$, $a_{12,10}=10$. The eigenvalues of the Laplacian matrix $L$ are $\lambda_1= \lambda_2 = \lambda_3 = 0$, $\lambda_4 = 1.101$, $\lambda_5 = 2.265$, $\lambda_6 = 3$, $ \lambda_{7,8} = 4 + 4i$, $\lambda_9 = 7$, $\lambda_{10} = 10.899$, $\lambda_{11} = 22$, $\lambda_{12} = 39.735$. It is observed that the Laplacian matrix has three zero eigenvalues that is consistent with the point that the underlying graph has three reach sets. Denote the right and the left eigenvectors associated with the zero eigenvalue of $L$  by $\gamma_1,\gamma_2,\gamma_3 \in \mathbb{R}^{12}$ and $\mu_1,\mu_2,\mu_3 \in \mathbb{R}^{12}$, respectively. In order to choose $\gamma_1$, we follow the properties~$1$ and $2$ in Lemma~\ref{lemma:eigenvectors_properties}, leading to $[\gamma_1]_1 = [\gamma_1]_2 = 1$ and $ [\gamma_1]_6 = [\gamma_1]_7 = 0$. The remaining entries of $\gamma_1$ are then obtained by solving $L \gamma_1 = \mathbf{0}_{12}$ that results in $\gamma_1 = [1~1~\mathbf{0}_{5}^T~\frac{1}{4}~\frac{1}{4}~\frac{1}{4}~\frac{1}{4}~\frac{1}{4}]^T$. It can be seen that $[\gamma_1]_i \in (0,1) $ for $i \in   \mathcal{C}_1$~(properties $3$ in Lemma~\ref{lemma:eigenvectors_properties}). Similarly, we have $\gamma_2 = [\mathbf{0}_{2}^T~1~1~1~\mathbf{0}_{2}^T~\frac{3}{4}~\frac{3}{4}~\frac{3}{4}~\frac{3}{4}~\frac{3}{4}]^T$ and $\gamma_3 =[\mathbf{0}_{5}^T~1~1~\mathbf{0}_{5}^T]^T$. It is now straightforward to check the property $4$ in Lemma~\ref{lemma:eigenvectors_properties}, i.e. $\gamma_1 + \gamma_2 + \gamma_3 = \mathbf{1}_{12}$. Turning into the construction of $\mu_1$, we use the second property in Lemma~\ref{lemma:eigenvectors_properties} that leads to $[\mu_1]_i = 0$ for $ i \in \mathcal{V} \backslash \mathcal{R}_1 = \{3,\dots,12\}$. The remaining entries of $\mu_1$ are obtained by solving $\mu_1^T L  = \mathbf{0}_{12}^T$ that results in $\mu_1 = [\frac{1}{3}~\frac{2}{3}~\mathbf{0}_{10}^T]^T$. It can be observed that $[\mu_1]_i \in (0,1) $ for $i \in   \mathcal{U}_1$~(properties $1$ in Lemma~\ref{lemma:eigenvectors_properties}), and $[\mu_1]_1 + [\mu_1]_2 = 1$~(properties $3$ in Lemma~\ref{lemma:eigenvectors_properties}). Following the same approach, we have $\mu_2=[\mathbf{0}_{2}^T~1~\mathbf{0}_{9}^T]^T$ and $\mu_2=[\mathbf{0}_{6}^T~1~\mathbf{0}_{5}^T]^T$.
\end{examcont}

\section{Problem Formulation and Results}
\label{sec:problem_results}

Motivated by the applications stated in Section~\ref{sec:introduction}, we investigate the effect of negative weights on the eigenvalues of the Laplacian matrix. Subsection~\ref{subsec:single_negative_weight} is mainly concerned with finding an upper bound on the magnitude of negative weights for added edges such that the eigenvalues of the Laplacian matrix have a specific property. We also present a discussion on the applicability of the notion of effective resistance for directed graphs. In Subsection~\ref{subsec:multiple_negative_weights}, we analyze the eigenvalues of Laplacian matrix for a graph with multiple negative weights whose magnitudes are infinitesimal.

\subsection{Adding an extra (un)directed negative edge to directed signed graphs}
\label{subsec:single_negative_weight}

In this subsection, we proceed to answering the following question:

\begin{question}
	\label{Question_1}Consider a signed graph~$\mathcal{G}_1(\mathcal{V},\mathcal{E}_1, \mathcal{W}_1)$ with the Laplacian matrix $L_1$. Construct a new graph $\mathcal{G}(\mathcal{V},\mathcal{E}, \mathcal{W}) = \mathcal{G}_1(\mathcal{V},\mathcal{E}_1, \mathcal{W}_1) \bigoplus \mathcal{G}_2(\mathcal{V},\mathcal{E}_2, \mathcal{W}_2)$,
	where $ \mathcal{E}_2= \{ (u,v), (v,u)\}$ and $\mathcal{W}_2(u,v) = - \delta_{uv} \leq 0$ and $ \mathcal{W}_2(v,u)= -\delta_{vu} \leq 0$. Denote $L$ the Laplacian matrix of $\mathcal{G}$.
	Find conditions on $\mathcal{G}_1$ and a bound on $\delta_{uv}$ and $\delta_{vu}$ such that
	\begin{equation}
		\label{eq:desired_eigenvalues_condition}
		0 = \lambda_1(L)  < \Re (\lambda_2(L)) \leq  \dots \leq \Re(\lambda_N(L)).
	\end{equation}
\end{question}

We first state the following definition and lemma that are used to prove the main results of this subsection.

\begin{definition}
	\label{def:r_delta}
	Consider a signed graph~$\mathcal{G}(\mathcal{V},\mathcal{E}, \mathcal{W})$ with the Laplacian matrix $L$. For given two arbitrary nodes $u$ and $v$, and variables $\delta_{uv},\delta_{vu},\omega \geq 0$, we define $r(\omega,\delta_{uv},\delta_{vu})$ as follows
	\begin{equation}
	\label{eq:r_delta_w}
	 r(\omega,\delta_{uv},\delta_{vu}) := (\mathbf{e}_u - \mathbf{e}_v)^T Q^T ( \bar{L}- j \omega I )^{-1} Q (\delta_{uv} \mathbf{e}_u - \delta_{vu} \mathbf{e}_v),
	\end{equation}
	where $\bar{L}=QLQ^T$ is the reduced Laplacian matrix.
\end{definition}

\begin{lemma}
	\label{lemma:eigenvalues_relation}
	Assume for a given directed graph $\mathcal{G}(\mathcal{V},\mathcal{E}, \mathcal{W})$ with Laplacian matrix $L$, $$\text{spec} \{ L \} = \{\underbrace{0,\dots,0}_{d}, \Lambda \},$$ where the set $\Lambda$ contains the non-zero eigenvalues of $L$. Then,
	\begin{enumerate}
		\item {$\text{spec} \{\bar{L} \} = \{\underbrace{0,\dots,0}_{d-1}, \Lambda \} $. Furthermore, if $d=1$, then $\bar{L}$ is invertible.
		}
		\item{Suppose further that all weights are non-negative.  $\bar{L}$ is invertible and $\Re \{ \lambda_i (\bar{L}) \} >0$ for $i=1,\dots,N-1$ if and only if $\mathcal{G}$ is connected~\footnote{For a graph with non-negative weights, the definition of connectivity, falls somewhere between the definitions of strong and weak connectivity~\cite{agaev2005spectra,young2015new}.}.
		}
	\end{enumerate}
\end{lemma}

\textit{Proof.}
	The proof of the first part can be found in~\cite[Lemma1]{young2010robustness}. Even though the lemma deals with directed graphs with non-negative weights, the same arguments are applicable to directed signed graphs, since $L \mathbf{1}_N = \mathbf{0}$. If $\mathcal{G}$ is connected with non-negative weights, then $\text{Spec}(L)=\{0,\lambda_2,\dots,\lambda_N\}$ with $\Re\{ \lambda_2 \} >0$~\cite{agaev2005spectra}. Taking into account this point along with the first part of the lemma proves the second part. \qed

We now present a necessity result  for the non-zero eigenvalues of Laplacian matrix to have positive real parts in the presence of a negative weight.

\begin{theorem}
	\label{theorem:one_negative_weight_undirected_necessity}
	Consider a signed graph~$\mathcal{G}_1(\mathcal{V},\mathcal{E}_1, \mathcal{W}_1)$ with the Laplacian matrix $L_1$. Assume that $L_1$ has only one zero eigenvalue and the rest of its eigenvalues have positive real parts.
	Construct a new graph $\mathcal{G}(\mathcal{V},\mathcal{E}, \mathcal{W}) = \mathcal{G}_1(\mathcal{V},\mathcal{E}_1, \mathcal{W}_1) \bigoplus \mathcal{G}_2(\mathcal{V},\mathcal{E}_2, \mathcal{W}_2)$,
	where $ \mathcal{E}_2= \{ (u,v), (v,u)\}$ and $\mathcal{W}_2(u,v) = -\delta_{uv} ,~\mathcal{W}_2(v,u)= -\delta_{vu}$ with $\delta_{uv} \geq 0,~\delta_{vu} \geq 0 $. Denote $L$ the Laplacian matrix of $\mathcal{G}$. If the eigenvalues of $L$ satisfy~(\ref{eq:desired_eigenvalues_condition}), then
	\begin{equation}
	\label{eq:cond_single_weight_general_necessity}
	\begin{aligned}
	r(0,\delta_{uv},\delta_{vu})  <1.
	\end{aligned}
	\end{equation}
\end{theorem}
\textit{Proof.}
	Define $\bar{L}= Q L Q^T$ and suppose that the eigenvalues of $L$ satisfy~(\ref{eq:desired_eigenvalues_condition}). We show the condition ~(\ref{eq:cond_single_weight_general_necessity}) holds. In view of the first part of Lemma~\ref{lemma:one_negative_weight_undirected} in~\ref{Appendix:a}, if $r(0,\delta_{uv},\delta_{vu})  > 1$, then $\det{ \left( \bar{L}_1^{-1} \bar{L} \right) } = \det{(\bar{L}_1^{-1})} \det{ (\bar{L})} < 0$. Hence $\bar{L}$ has at least one non-positive eigenvalue. If $r(0,\delta_{uv},\delta_{vu})  = 1$, then $\det{ \left( \bar{L}_1^{-1} \bar{L} \right) } = 0$. Thus, $\bar{L}$ has at least one zero eigenvalue. This means that the condition~(\ref{eq:cond_single_weight_general_necessity}) is necessary for the eigenvalues of $L$ satisfy~(\ref{eq:desired_eigenvalues_condition}). \qed

We now present a sufficiency theorem for the non-zero eigenvalues of Laplacian matrix to have positive real parts in the presence of a negative weight.

\begin{theorem}
	\label{theorem:one_negative_weight_undirected_sufficiency}
	Consider a signed graph~$\mathcal{G}_1(\mathcal{V},\mathcal{E}_1, \mathcal{W}_1)$ with the Laplacian matrix $L_1$. Assume that $L_1$ has only one zero eigenvalue and the rest of its eigenvalues have positive real parts.
	Construct a new graph $\mathcal{G}(\mathcal{V},\mathcal{E}, \mathcal{W}) = \mathcal{G}_1(\mathcal{V},\mathcal{E}_1, \mathcal{W}_1) \bigoplus \mathcal{G}_2(\mathcal{V},\mathcal{E}_2, \mathcal{W}_2)$,
	where $ \mathcal{E}_2= \{ (u,v), (v,u)\}$ and $\mathcal{W}_2(u,v) = -\delta q_{uv} \leq 0 ,~\mathcal{W}_2(v,u)= -\delta q_{vu} \leq 0$ with $\delta > 0 $ and given $q_{uv}, q_{vu} \geq 0$. Denote $L$ the Laplacian matrix of $\mathcal{G}$. Let $\delta^*$ be obtained by,
	\begin{equation}
	\label{eq:cond_single_weight_general_sufficiency_final}
	\begin{aligned}
	\min_{\delta_1 \in \mathbb{R}_{> 0}, \omega \in \mathbb{R}_{\geq 0}} &~~~~~ \delta_1 \\
	\text{subject to} &~~~~~ r(\omega,\delta_1 q_{uv},\delta_1 q_{vu}) = 1.
	\end{aligned}
	\end{equation}
	Then, the eigenvalues of $L$ satisfy~(\ref{eq:desired_eigenvalues_condition}) for all $ \delta \in [0,\delta^*)$.
\end{theorem}

\textit{Proof.} Denote $L_2$ the Laplacian matrix of $\mathcal{G}_2$. Since $\mathcal{E}_2= \{ (u,v), (v,u)\}$, $L_2$ can be expressed as $L_2 = - (\delta_{uv} \mathbf{e}_u - \delta_{vu} \mathbf{e}_v) (\mathbf{e}_u - \mathbf{e}_v)^T$ with $\delta_{uv} = \delta q_{uv}$ and $\delta_{vu} = \delta q_{vu}$. Under conditions of the theorem, $L$ can be written as $L=L_1 + L_2$,  leading to the following expression for $\bar{L}$,
\begin{equation}
\label{eq:Lbar_new_representation_new}
\bar{L} =\underbrace{Q L_1 Q^T}_{\bar{L}_1} - Q (\delta_{uv} \mathbf{e}_u - \delta_{vu} \mathbf{e}_v)  (\mathbf{e}_u - \mathbf{e}_v)^T Q^T.
\end{equation}
Since $L_1$ has only one zero eigenvalue and the rest of its eigenvalues have positive real parts, all eigenvalues of $\bar{L}_1 - j \omega I$ have positive real parts  according to Property $1$ in Lemma~\ref{lemma:eigenvalues_relation}.

We now prove the theorem for the case with $q_{uv}, q_{vu} > 0$. The proof of theorem for the cases with $q_{uv} > 0,~q_{vu}=0$ or $q_{uv} =0,~q_{vu} > 0$  follows the same lines. Since $q_{uv}, q_{vu} > 0$, we have $\delta_{uv} = \delta q_{uv}$ and $\delta_{vu} = \delta q_{vu}$.  Note that  for $\delta = 0$, $\bar{L} = \bar{L}_1$ and for sufficiently large $\delta$, i.e. $\delta > \sum\limits_{i=1}^{N} [L_1]_{ii}$, the sum of the diagonal entries of $L$ becomes negative. This leads $L$ to have at least one eigenvalue with a negative real part, and consequently, in view of Property 1 in Lemma~\ref{lemma:eigenvalues_relation}, $\bar{L}$ has that eigenvalue with negative real part. Hence, the continuity of eigenvalues of $\bar{L}$ with respect to $\delta$ states that
\begin{equation}
\label{eq:delta_star_new}
\exists~\delta^*>0~\text{such that}~\forall~\delta \in [0,\delta^*),~ \Re \{ \lambda_i ( \bar{L} ) \}_{i=1}^{N-1} > 0,
\end{equation}
and also with $\delta = \delta^*$, $\Re \{ \lambda_i ( \bar{L} ) \} = 0$ for some $i$. This means there exists at least one $\omega \in \mathbb{R}_{\geq 0}$ such that $\lambda_i ( \bar{L} ) = j \omega  $ for $\delta = \delta^*$, or equivalently $\det{ \left( \bar{L}- j \omega I \right) } = 0$. Using the expression of $\bar{L}$ in~(\ref{eq:Lbar_new_representation_new}) and taking into account that $\bar{L}_1 - j \omega I$ is invertible, we have
\begin{equation}
\label{eq:determinant_jw}
\begin{aligned}
\det{ \left( \bar{L}- j \omega I \right) }\ &=  \det{ \left( \bar{L}_1- j \omega I - Q (\delta_{uv} \mathbf{e}_u - \delta_{vu} \mathbf{e}_v) (\mathbf{e}_u - \mathbf{e}_v)^T Q^T \right) } \\
&= \det{ ( \bar{L}_1- j \omega I )} ~ { \left( 1-  (\mathbf{e}_u - \mathbf{e}_v)^T Q^T ( \bar{L}_1- j \omega I )^{-1} Q (\delta_{uv} \mathbf{e}_u - \delta_{vu} \mathbf{e}_v) \right)}  \\
&= \det{ ( \bar{L}_1- j \omega I )} ~ { \left( 1- \delta^* (\mathbf{e}_u - \mathbf{e}_v)^T Q^T ( \bar{L}_1- j \omega I )^{-1} Q ( q_{uv}\mathbf{e}_u - q_{vu} \mathbf{e}_v) \right)}\\
&= 0,
\end{aligned}
\end{equation}
where the last two equalities are obtained by applying Lemma~\ref{lemma:rank_one_updates} and taking into account  $\delta_{uv} = \delta q_{uv}$ and $\delta_{vu} = \delta q_{vu}$. Since $\bar{L}_1 - j \omega I$ is invertible, we observe that $ \det{ \left( \bar{L}- j \omega I \right) } = 0$ if and only if $\delta^* (\mathbf{e}_u - \mathbf{e}_v)^T Q^T ( \bar{L}_1- j \omega I )^{-1} Q ( q_{uv}\mathbf{e}_u - q_{vu}\mathbf{e}_v) = 1$. By solving~(\ref{eq:cond_single_weight_general_sufficiency_final}), we find the minimum value of $\delta^*$ such that the eigenvalues of $L$ satisfy~(\ref{eq:desired_eigenvalues_condition}) for all $ \delta \in (0,\delta^*)$. The optimization problem~(\ref{eq:cond_single_weight_general_sufficiency_final}) always has a solution since  the continuity argument above guarantees that  the eigenvalues of $\bar{L}$ cross the imaginary axis for some $\delta \in (0,\delta_1 ]$.  This completes the proof of the case $q_{uv}, q_{vu} > 0$.\qed
%

\begin{remark}
	Theorem~\ref{theorem:one_negative_weight_undirected_sufficiency} includes three different cases that correspond to perturbing the edge $(u,v)$~($q_{uv}> 0,q_{vu}=0$), the edge~$(v,u)$~($q_{uv}=0,q_{vu}> 0$), or both edges~$(u,v)$ and~$(v,u)$~($q_{uv}>0,q_{vu}>0$). The variables $q_{uv},q_{vu}$ allow to incorporate all these cases in only one optimization problem as stated in~(\ref{eq:cond_single_weight_general_sufficiency_final}). 
\end{remark}

The sufficiency condition~(\ref{eq:cond_single_weight_general_sufficiency_final}) in~Theorem~\ref{theorem:one_negative_weight_undirected_sufficiency} correspond to the minimum value of $\delta^*$ such that the eigenvalues of $L$ satisfy~(\ref{eq:desired_eigenvalues_condition}) for all $ \delta \in (0,\delta^*)$. One of the following statements holds for $\delta^*$ obtained from solving~(\ref{eq:cond_single_weight_general_sufficiency_final}):

\begin{enumerate}
	\item {$\delta^*$ is obtained with $\omega_1,\dots,\omega_k$ where $\omega_i \neq 0$ for $i=1,\dots,k$ and $k$ is the number of solutions for~(\ref{eq:cond_single_weight_general_sufficiency_final}) which is finite since the constraint is a non trivial polynomial equation in $\omega$ by definition of $r(\omega,\delta_{uv},\delta_{vu})$;
	}
	\item {or, there exists at least one zero $\omega_i$.
	}
\end{enumerate}

In the second case, the matrix $\bar{L}$ has at least one zero eigenvalue with $\delta = \delta^*$~(in the view of~(\ref{eq:delta_star_new}) ) which means that $\det{ \left( \bar{L}_1^{-1} \bar{L} \right) } = 0$ for $\delta = \delta^*$.  In this case,  the condition~(\ref{eq:cond_single_weight_general_necessity}) becomes necessary and sufficient for the non-zero eigenvalues of Laplacian matrix have positive real parts in the presence of a negative weight according to the following theorem.  

\begin{theorem}
	\label{theorem:one_negative_weight_undirected_corrected}
	Consider a signed graph~$\mathcal{G}_1(\mathcal{V},\mathcal{E}_1, \mathcal{W}_1)$ with the Laplacian matrix $L_1$. Assume that $L_1$ has only one zero eigenvalue and the rest of its eigenvalues have positive real parts.
	Construct a new graph $\mathcal{G}(\mathcal{V},\mathcal{E}, \mathcal{W}) = \mathcal{G}_1(\mathcal{V},\mathcal{E}_1, \mathcal{W}_1) \bigoplus \mathcal{G}_2(\mathcal{V},\mathcal{E}_2, \mathcal{W}_2)$,
	where $ \mathcal{E}_2= \{ (u,v), (v,u)\}$ and $\mathcal{W}_2(u,v) = -\delta q_{uv} \leq 0 ,~\mathcal{W}_2(v,u)= -\delta q_{vu} \leq 0$ with $\delta > 0 $ and given $q_{uv}, q_{vu} \geq 0 $. Denote $L$ the Laplacian matrix of $\mathcal{G}$.
  Let $\delta^*$ be obtained from~(\ref{eq:cond_single_weight_general_sufficiency_final}) with $\omega_1,\dots,\omega_k$ being the roots of the equality constraint. Assume further that there exists at least one zero $\omega_i$. The eigenvalues of $L$ satisfy~(\ref{eq:desired_eigenvalues_condition}) if and only if
		\begin{equation}
		\label{eq:cond_single_weight_general_final}
		\begin{aligned}
		\underbrace{\delta \left( \mathbf{e}_u - \mathbf{e}_v \right)^T    Q^T  \bar{L}_1^{-1} Q \left( q_{uv}\mathbf{e}_u - q_{vu} \mathbf{e}_v \right)}_{r(0,\delta q_{uv},\delta q_{vu})} <1.
		\end{aligned}
		\end{equation}
		Furthermore, $\delta^* \left( \mathbf{e}_u - \mathbf{e}_v \right)^T    Q^T  \bar{L}_1^{-1} Q \left( q_{uv} \mathbf{e}_u - q_{vu} \mathbf{e}_v \right)  = 1$ and $ \left( \mathbf{e}_u - \mathbf{e}_v \right)^T    Q^T  \bar{L}_1^{-1} Q \left( q_{uv} \mathbf{e}_u - q_{vu} \mathbf{e}_v \right) >0  $.
\end{theorem}

\textit{Proof.}	The necessity part results from Theorem~\ref{theorem:one_negative_weight_undirected_necessity}, i.e. if the eigenvalues of $L$ satisfy~(\ref{eq:desired_eigenvalues_condition}), then~(\ref{eq:cond_single_weight_general_final}) holds. We now show that the condition~(\ref{eq:cond_single_weight_general_final}) is sufficient for the case with  $ q_{uv}, q_{vu} > 0$. The cases with $q_{uv} >0,~q_{vu} = 0$ or $q_{uv} = 0,~q_{vu} >0$ can be proven by using the same arguments.
\newline
In the proof of Theorem~\ref{theorem:one_negative_weight_undirected_sufficiency}, it is observed from the continuity of eigenvalues of $\bar{L}$ with respect to $\delta$ that
\begin{equation}
\label{eq:delta_star}
\exists~\delta^*>0~\text{such that}~\forall~\delta \in [0, \delta^*),~ \Re \{ \lambda_i ( \bar{L} ) \}_{i=1}^{N-1} > 0,
\end{equation}
and also with $\delta = \delta^*$, $\Re \{ \lambda_i ( \bar{L} ) \} = 0$ for some $i$. Since $\delta^*$ is obtained from~(\ref{eq:cond_single_weight_general_sufficiency_final}) with at least one zero $\omega_i$, then there exists at least one zero eigenvalue~$  \lambda_i ( \bar{L} )  = 0$ for $\delta = \delta^*$, meaning that $\det{ \left( \bar{L}_1^{-1} \bar{L} \right) } = 0$. From the definition of $r(0,\delta q_{uv},\delta q_{vu})$ in~(\ref{eq:cond_single_weight_general_final}), we have
\begin{equation}
\label{eq:r_delta}
r_{\delta} := r(0,\delta q_{uv},\delta q_{vu})= \alpha \delta,
\end{equation}
where $\alpha = \left( \mathbf{e}_u - \mathbf{e}_v \right)^T    Q^T  \bar{L}_1^{-1} Q \left( q_{uv}\mathbf{e}_u - q_{vu} \mathbf{e}_v \right)$.
\newline
Let assume $\alpha > 0$. Using this assumption, (\ref{eq:delta_star}) can be rewritten as
\begin{equation}
\label{eq:r_star}
\exists~r^*> 0~\text{such that}~\forall~r_{\delta} \in [0, r^*),~ \Re \{ \lambda_i ( \bar{L} ) \}_{i=1}^{N-1} > 0.
\end{equation}
We now need to show $r^* =1$. If $ r^* < 1$, according to~(\ref{eq:r_star}), the real part of at least one eigenvalue of $\bar{L} $ is negative for $r_{\delta} > r^*$ and the real part becomes zero at $r_{\delta} = r^*$; however, from Lemma~\ref{lemma:one_negative_weight_undirected} in~\ref{Appendix:a}, it is observed that $\bar{L}$ can have zero eigenvalue(s) if and only if $ r_{\delta} = 1$. This contradicts the assumption $ r^* < 1$. Hence, $r^* = 1$
\newline
To complete the proof, we should show the assumption $\alpha > 0 $. To do so, we show $\alpha$ cannot be zero or strictly negative by considering two cases.

Case 1: To obtain a contradiction, assume $\alpha = 0$. Then $r_{\delta} = 0$, and  according to~(\ref{eq:spectrum_L_1bar_L_bar}), $\det{ \left( \bar{L}_1^{-1} \bar{L} \right) }$ cannot be zero which means $\bar{L}$ has no zero eigenvalue. This means that $\delta^* = +\infty$ or equivalently all eigenvalues of $\bar{L}$ have positive real parts. However, this is not true, since for sufficiently large $\delta$, i.e. $\delta > \sum\limits_{i=1}^{N} [L_1]_{ii}$, the sum of the diagonal entries of $L$ becomes negative. This leads $L$ to have at least one eigenvalue with a negative real part. Hence, in view of Property 1 in Lemma~\ref{lemma:eigenvalues_relation}, $\bar{L}$ has the same eigenvalue with negative real part. This is in contradiction with $\delta^* = + \infty$. Therefore, $\alpha \neq 0$.

Case 2: To obtain a contradiction, assume $\alpha < 0$. Thus (\ref{eq:delta_star}) and (\ref{eq:r_delta}) imply
\begin{equation}
\label{eq:r_star_1}
\exists~r^*>0~\text{such that}~\forall~r_{\delta} \in (-r^*, 0],~ \Re \{ \lambda_i ( \bar{L} ) \}_{i=1}^{N-1} > 0.
\end{equation}
Furthermore, letting $r_{\delta} = -r^*$ yields  $\Re \{ \lambda_i ( \bar{L} ) \} = 0$ for some $i$  as $r^*$ is taken as the maximum value satisfying~(\ref{eq:r_star_1}). This means that $\det{ \left( \bar{L}_1^{-1} \bar{L} \right) } = 0 $ for $r_{\delta} = -r^* < 0$. However, according to ~(\ref{eq:spectrum_L_1bar_L_bar}), $\det{ \left( \bar{L}_1^{-1} \bar{L} \right) } = 0 $ if and only if $r_\delta=1$. This contradicts $r_{\delta} < 0$. Hence, $\alpha $ cannot be negative. This completes the proof. \qed

The key step to apply Theorems~\ref{theorem:one_negative_weight_undirected_sufficiency} and~\ref{theorem:one_negative_weight_undirected_corrected} is to numerically solve the optimization problem~(\ref{eq:cond_single_weight_general_sufficiency_final}), which has at least one solution as explained in the proof of Theorem~\ref{theorem:one_negative_weight_undirected_sufficiency}. The feasible set of this optimization problem can be interpreted in terms of the Nyquist plots of the following system~$\Sigma_{uv}$,
\begin{equation}
\label{eq:dynamic_systems}
\begin{aligned}
&\Sigma_{uv} := \left\{
\begin{array}{ll}
\frac{d}{dt}x_{uv} &= \bar{L}_1 x_{uv} + Q ( q_{uv} \mathbf{e}_u - q_{vu} \mathbf{e}_v) u_{uv}\\
y_{uv} &= (\mathbf{e}_u - \mathbf{e}_v)^T Q^T x_{uv}
\end{array}
\right.\\
\end{aligned}
\end{equation}
where $x_{uv} \in \mathbb{R}^{N-1}$, $u_{uv} \in \mathbb{R}$ and $y_{uv} \in \mathbb{R}$  denote the state, the input and the output of the system $\Sigma_{uv}$.
Consider the optimization problem~(\ref{eq:cond_single_weight_general_sufficiency_final}) and the system $\Sigma_{uv}$ in a negative feedback structure with a proportional controller $\delta_1$. The optimization problem~(\ref{eq:cond_single_weight_general_sufficiency_final}) yields the minimum value of the gain for which the Nyqiust diagram of $ \delta_1 G_{uv}(s)$ crosses the critical point $-1$ where $  G_{uv}(s)$ is the transfer function from $u_{uv}$ to $y_{uv}$. Hence, to attain $\delta^*$, one can plot the Nyquist diagram of $  \Sigma_{uv}$ and find frequencies $\omega_i,~i=1,\dots,k$, at which it crosses the real axis. Then, $\delta^* = \frac{1}{|G_{uv}(j\omega^*)|}$ where $G_{uv}(j\omega^*) \leq G_{uv}(j\omega_i)$ for $i=1,\dots,k$. If $\omega^* = 0$, then we can use the results of  Theorem~\ref{theorem:one_negative_weight_undirected_corrected}; otherwise we should use the results of  Theorem~\ref{theorem:one_negative_weight_undirected_sufficiency}. In the following examples, we illustrate how to apply the results of Theorems~\ref{theorem:one_negative_weight_undirected_sufficiency} and~\ref{theorem:one_negative_weight_undirected_corrected}.

\begin{figure}[H]
	\centering
	\includegraphics[scale=0.4]{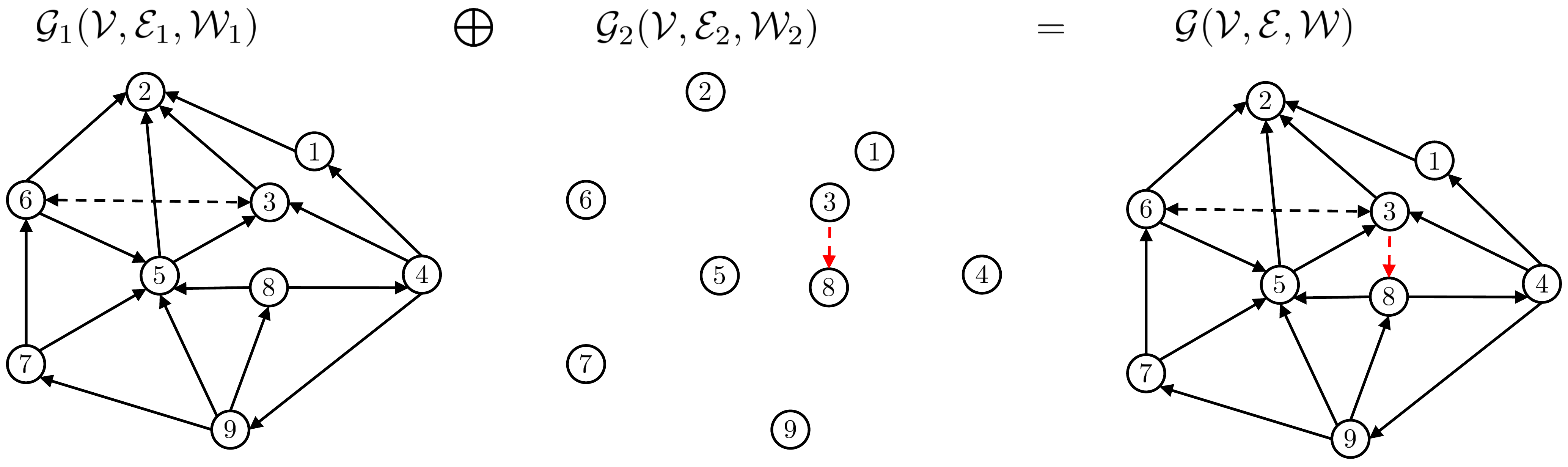}
	\caption{A directed graph studied in Example~1. The solid arrows represent the positive weight edges while the dashed arrows show the edges with negative weights. All positive and negative weights are set equal to $2$ and $-1$, respectively.}
	\label{fig:fig2}
\end{figure}

\begin{figure}[H]
	\centering
	\includegraphics[scale=0.4]{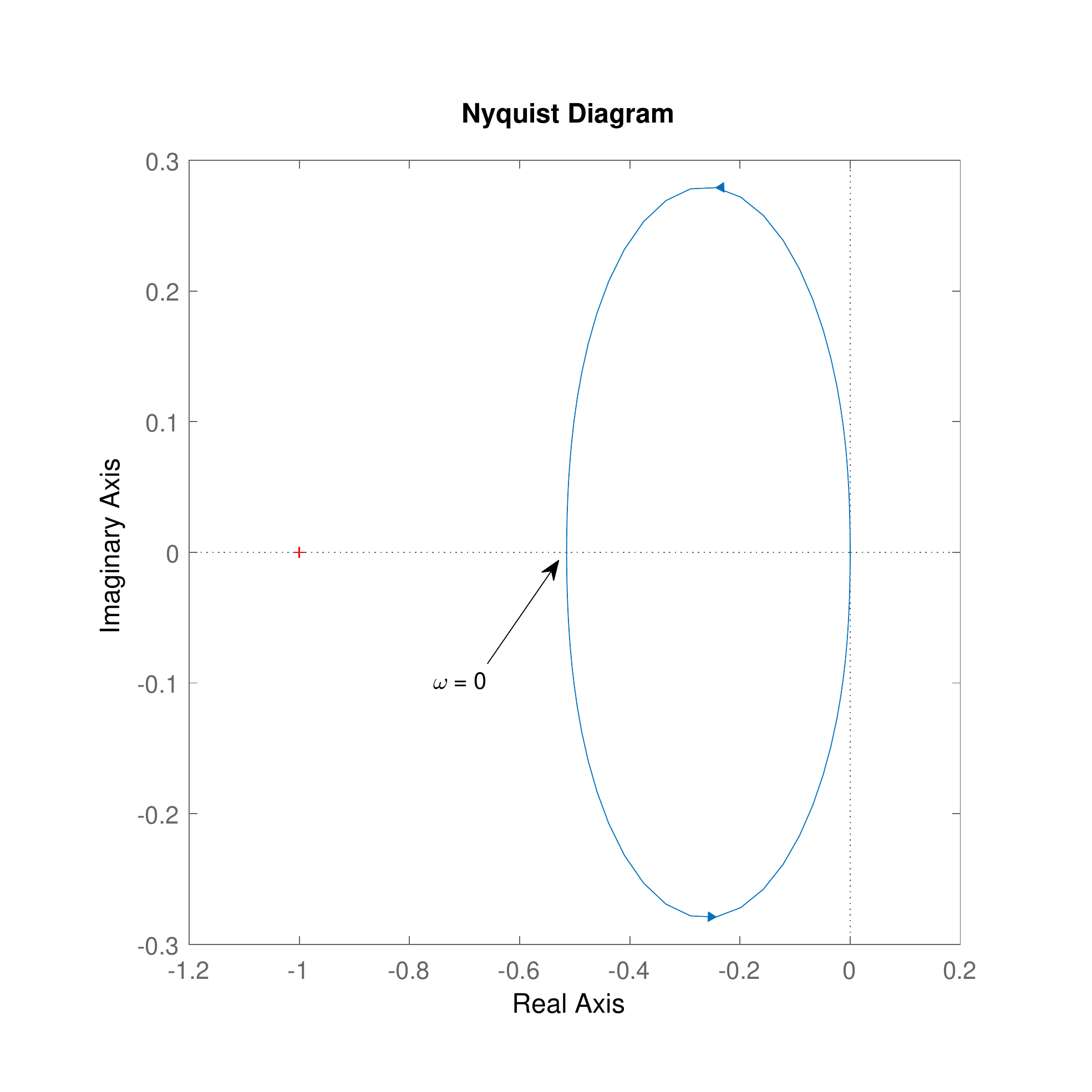}
	\caption{The Nyquist diagram in Example~$1$.}
	\label{fig:fig2_Nyquist}
\end{figure}

\begin{example}
	Consider the graphs $\mathcal{G}_1(\mathcal{V},\mathcal{E}_1, \mathcal{W}_1)$, $\mathcal{G}_2(\mathcal{V},\mathcal{E}_2, \mathcal{W}_2)$ and~$\mathcal{G}(\mathcal{V},\mathcal{E}, \mathcal{W})$ depicted in Figure~\ref{fig:fig2} with the Laplacian matrices $L_1$, $L_2$, $L$ respectively. In this figure, the solid arrows represent the positive weights which are set equal to $2$, and  the dashed arrows represent the negative weights with $a_{36} = a_{63} = a_{38} = -1$. The graph $\mathcal{G}$ is constructed from the graphs $\mathcal{G}_1$ and $\mathcal{G}_2$ by adding negative weights between two pairs of nodes, i.e. $\mathcal{G} = \mathcal{G}_1 \oplus \mathcal{G}_2$. We now add the negative directed edge $a_{38} = -\delta$. The Nyquist diagram plotted in Figure~\ref{fig:fig2_Nyquist} crosses the real axis at $\omega_1 = 0$ and $\omega_2 = \infty$. Furthermore, the magnitude of the Nyquist diagram at $\omega_1$ is smaller than of that at $\omega_2$. Hence, in view of Theorem~\ref{theorem:one_negative_weight_undirected_corrected}, the condition~(\ref{eq:cond_single_weight_general_final}) is necessary and sufficient condition, leading to $\delta^*= 1.94285 $. This means that the eigenvalues of $L$  satisfy~(\ref{eq:desired_eigenvalues_condition}) if and only if $\delta < 1.94285 $.
\end{example}

\begin{figure}[H]
	\centering
	\includegraphics[scale=0.4]{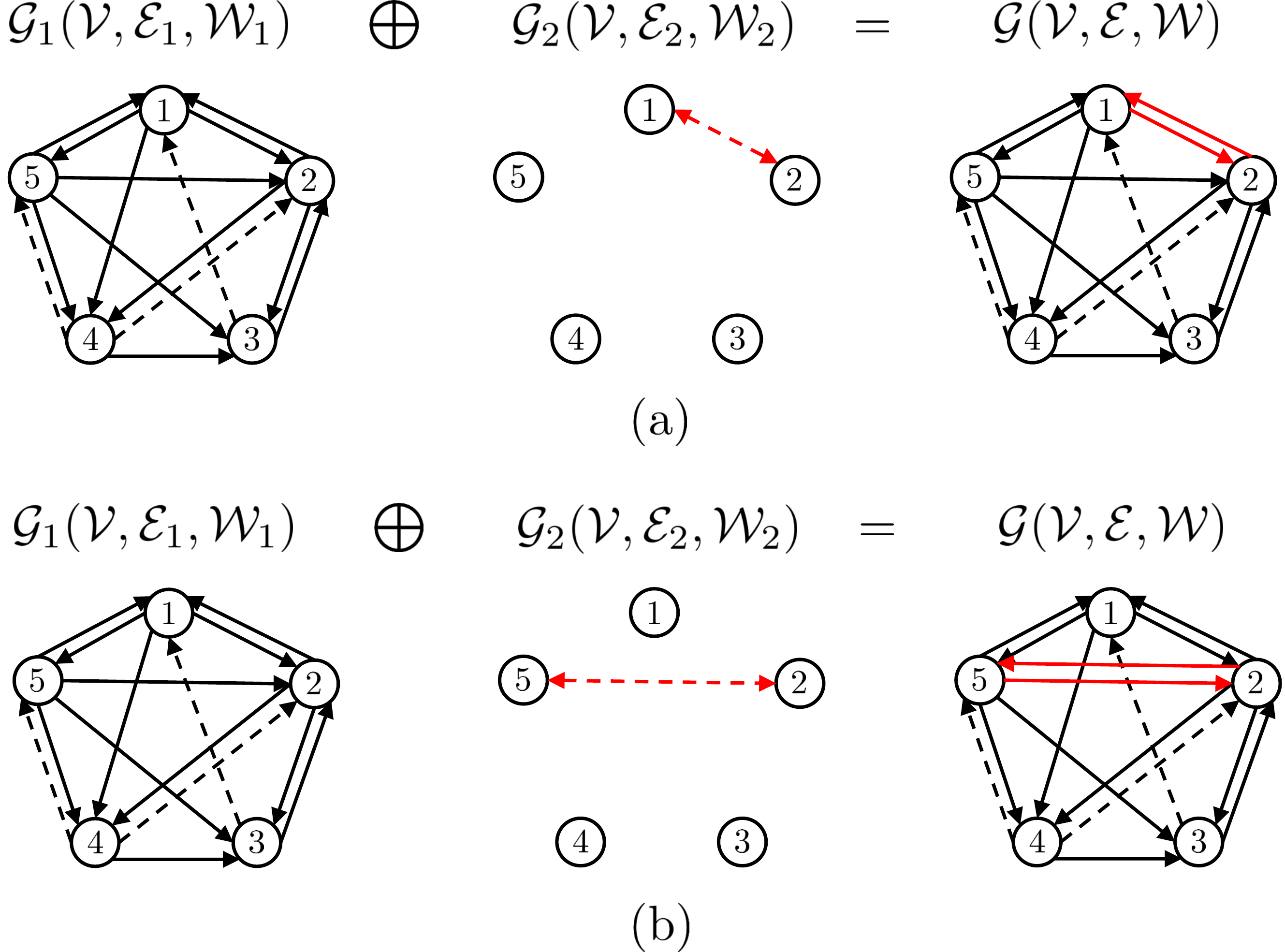}
	\caption{A directed graph studied in Example~2. The solid arrows represent the positive weight edges while the dashed arrows show the edges with negative weights. All positive and negative weights are set equal to $2$ and $-1$, respectively.}
	\label{fig:fig2_1}
\end{figure}

\begin{figure}[H]
	\centering
	\centerline{%
		\includegraphics[width=0.5\textwidth]{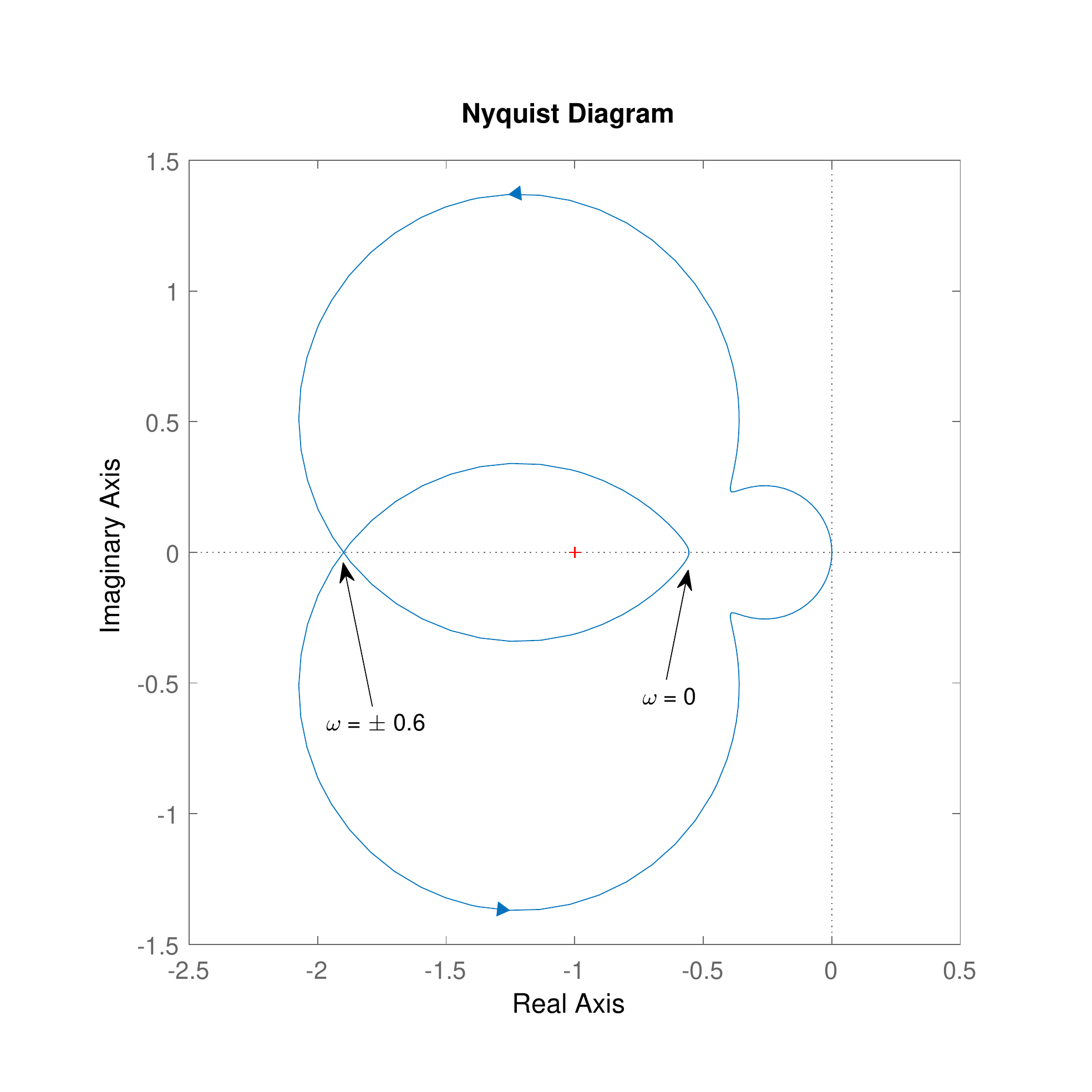}%
		\includegraphics[width=0.5\textwidth]{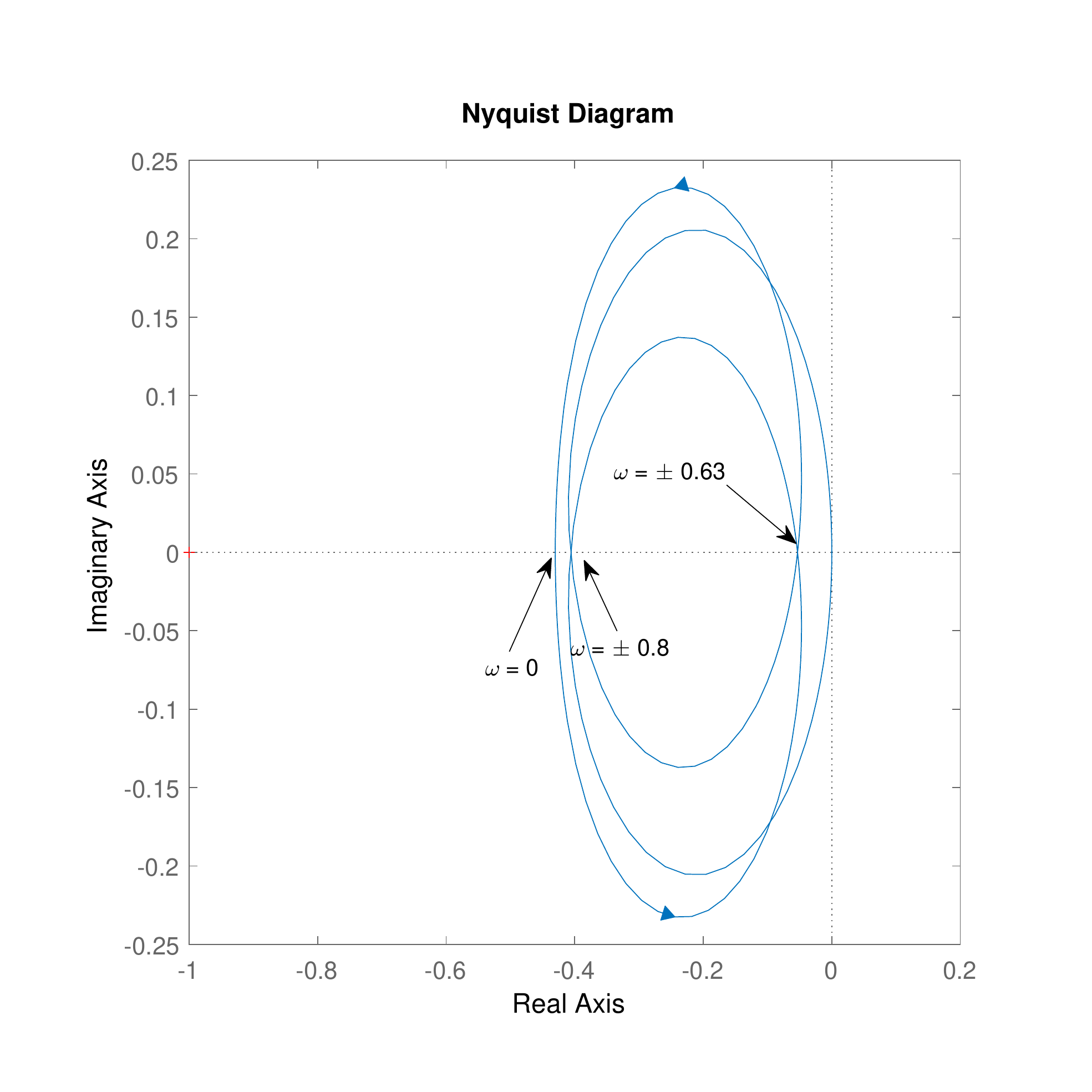}%
	}%
	\caption{The Nyquist diagrams in Example~$2$. The left Nyquist diagram corresponds to the first scenario with the graph in Figure~\ref{fig:fig2_1}(a), while the right Nyquist diagram corresponds to the second scenario with the graph in Figure~\ref{fig:fig2_1}(b). }
	\label{fig:fig2_1_Nyquist}
\end{figure}

\begin{example}
	Consider the graphs $\mathcal{G}_1(\mathcal{V},\mathcal{E}_1, \mathcal{W}_1)$, $\mathcal{G}_2(\mathcal{V},\mathcal{E}_2, \mathcal{W}_2)$ and~$\mathcal{G}(\mathcal{V},\mathcal{E}, \mathcal{W})$ depicted in Figure~\ref{fig:fig2_1} with the Laplacian matrices $L_1$, $L_2$, $L$ respectively. In this figure, the solid black arrows, the dashed red arrow and the solid read arrows represent the positive weights,  the negative weight and the perturbed edges, respectively.  We assign the following weights to the edges of the graph $\mathcal{G}_1$: $a_{12}=1$, $a_{14}=1$, $a_{15}=1$, $a_{21}=1$, $a_{23}=1$, $a_{24}=1$, $a_{31}=-1$, $a_{32}=1$, $a_{35}=-0.8$, $a_{42}=-0.3$, $a_{43}=1.5$, $a_{45}=-2$, $a_{51}=2$, $a_{52}=1$, $a_{53}=2$, $a_{54}=1$.
	\newline
	In the first scenario, the graph $\mathcal{G}$ is constructed from the graphs $\mathcal{G}_1$ and $\mathcal{G}_2$ by perturbing positive weights between two pairs of nodes $1$ and $2$. According to the left Nyquist diagram plotted in Figure~\ref{fig:fig2_1_Nyquist}, the Nyquist diagram crosses the real axis at $\omega_1 = 0$, $\omega_{2,3} = \pm 0.6$ and $\omega_{4,5} = \pm \infty$. Furthermore, the magnitude of the Nyquist diagram at $\omega_{2,3}$ is smaller than others. Hence, in view of Theorem~\ref{theorem:one_negative_weight_undirected_corrected}, the condition~(\ref{eq:cond_single_weight_general_final}) is sufficient, leading to $\delta^*= 0.52 $. This means that if $\delta < 0.52 $ the eigenvalues of $L$  meet~(\ref{eq:desired_eigenvalues_condition}). On the other hand, the necessary condition~(\ref{eq:cond_single_weight_general_necessity}) holds for $\delta < 1.8$.
	\newline
	In the second scenario, the graph $\mathcal{G}$ is constructed from the graphs $\mathcal{G}_1$ and $\mathcal{G}_2$ by perturbing positive weights between two pairs of nodes $2$ and $5$. According to the right Nyquist diagram plotted in Figure~\ref{fig:fig2_1_Nyquist}, the Nyquist diagram crosses the real axis at $\omega_1 = 0$, $\omega_{2,3} = \pm 0.63$, $\omega_{4,5} = \pm 0.8$ and $\omega_{6,7} = \pm \infty$. Furthermore, the magnitude of the Nyquist diagram at $\omega_1$ is smaller than others. Hence, in view of Theorem~\ref{theorem:one_negative_weight_undirected_corrected}, the condition~(\ref{eq:cond_single_weight_general_final}) is necessary and sufficient condition, leading to $\delta^*= 2.3239 $. This means that the eigenvalues of $L$  meet~(\ref{eq:desired_eigenvalues_condition}) if and only if $\delta < 2.3239 $.
\end{example}

\begin{remark}
	\label{remark:possitiveness_of_multiplications_correction}
	The results of this subsection are different from those in~\cite{mukherjee2016consensus} where sufficient conditions for the  upper bound $\delta$ has been derived via Nyquist stability criteria. First, we provide the necessary result which considers the case in which both edges between any arbitrary pairs of nodes are perturbed with  negative weights. Secondly, our sufficiency result is more general than the main result of~\cite{mukherjee2016consensus} since we also allow  perturbing two edges between two nodes with the same negative weight.  Our results cover a more general set of graphs as a graph with multiple negative edges might satisfy the assumption of the theorem, while~\cite[Theorem 1]{mukherjee2016consensus} only is applied to graphs with no negative edges. Even though the results of this paper are interpreted via Nyquist criteria, the definitions of systems are different from~\cite{mukherjee2016consensus}. Finally, we highlight the case where the condition becomes necessary and sufficient.
\end{remark}

In all theorems above, it is assumed that the original graph~$\mathcal{G}_1$ satisfy~(\ref{eq:desired_eigenvalues_condition}). If  $\mathcal{G}_1$ is connected and all of its weight are non-negative, the second part of Lemma~\ref{lemma:eigenvalues_relation} provides the necessary and sufficient conditions for $\mathcal{G}_1$ to satisfy~(\ref{eq:desired_eigenvalues_condition}). In this case, the following corollaries are obtained.

\begin{corollary}
	\label{corollary:one_negative_weight_undirected}
	Consider a signed graph~$\mathcal{G}(\mathcal{V},\mathcal{E}, \mathcal{W})$ with the Laplacian matrix $L$ where it is decomposed into two subgraphs
	$\mathcal{G} (\mathcal{V},\mathcal{E}^+,\mathcal{W}^+)$ and $\mathcal{G} (\mathcal{V},\mathcal{E}^-,\mathcal{W}^-)$ with the corresponding Laplacian matrices $L^+$ and $L^{-}$, respectively.
	Assume $\mathcal{E}^{-}= \{ (u,v), (v,u)\}$ with $\mathcal{W}^{-}(u,v) = - \delta q_{uv} \leq 0$ and $\mathcal{W}^{-}(v,u) = - \delta q_{vu} \leq 0 $ with $\delta > 0 $ and given $q_{uv}, q_{vu} \geq 0$. Assume also $\mathcal{G} (\mathcal{V},\mathcal{E}^+,\mathcal{W}^+)$ is connected, and define  $\bar{L}^+ = Q L^+ Q^T$. Let $\delta^*$ be obtained by,
	\begin{equation}
	\label{eq:cond_single_weight_general_connected_postive_sufficiency_final}
	\begin{aligned}
	\min_{ \omega \in \mathbb{R}_{\geq 0}}&~~~~~\delta_1 \\
	\text{subject to} &~~~~~\delta_1 (\mathbf{e}_u - \mathbf{e}_v)^T Q^T ( \bar{L}^+ - j \omega I )^{-1} Q ( q_{uv} \mathbf{e}_u - q_{vu} \mathbf{e}_v) = 1,~~~\delta_1 >0.
	\end{aligned}
	\end{equation}
	Then, the eigenvalues of $L$ satisfy~(\ref{eq:desired_eigenvalues_condition}) for all $ \delta \in [0,\delta^*)$.
	%
\end{corollary}
\begin{corollary}
	\label{corollary:one_negative_weight_undirected_1}
	Consider a signed graph~$\mathcal{G}(\mathcal{V},\mathcal{E}, \mathcal{W})$ with the Laplacian matrix $L$ where it is decomposed into two subgraphs
	$\mathcal{G} (\mathcal{V},\mathcal{E}^+,\mathcal{W}^+)$ and $\mathcal{G} (\mathcal{V},\mathcal{E}^-,\mathcal{W}^-)$ with the corresponding Laplacian matrices $L^+$ and $L^{-}$, respectively.
	Assume $\mathcal{E}^{-}= \{ (u,v), (v,u)\}$ with $\mathcal{W}^{-}(u,v) = - \delta q_{uv} \leq 0$ and $\mathcal{W}^{-}(v,u) = - \delta q_{vu} \leq 0 $ with $\delta > 0 $ and given $q_{uv}, q_{vu} \geq 0$. Assume also $\mathcal{G} (\mathcal{V},\mathcal{E}^+,\mathcal{W}^+)$ is connected, and define  $\bar{L}^+ = Q L^+ Q^T$. Let $\delta^*$ be obtained from~(\ref{eq:cond_single_weight_general_sufficiency_final}) with $\omega_1,\dots,\omega_k$ for $i=1,\dots,k$. Assume further that there exists at least one zero $\omega_i$. The eigenvalues of $L$ satisfy~(\ref{eq:desired_eigenvalues_condition}) if and only if
	\begin{equation}
	\label{eq:cond_single_weight_general_connected_postive_final}
	\begin{aligned}
	\delta \underbrace{ \left( \mathbf{e}_u - \mathbf{e}_v \right)^T    Q^T  \left(\bar{L}^+\right)^{-1} Q \left( q_{uv} \mathbf{e}_u - q_{vu} \mathbf{e}_v \right)}_{r(0, q_{uv}, q_{vu})} <1.
	\end{aligned}
	\end{equation}
	Furthermore, $\delta^* \left( \mathbf{e}_u - \mathbf{e}_v \right)^T    Q^T  \left(\bar{L}^+\right)^{-1} Q \left( \mathbf{e}_u - \mathbf{e}_v \right)  = 1$ and $ \left( \mathbf{e}_u - \mathbf{e}_v \right)^T    Q^T  \left(\bar{L}^+\right)^{-1} Q \left( \mathbf{e}_u - \mathbf{e}_v \right) >0  $.

\end{corollary}

\begin{figure}[h!]
	\centering
	\begin{subfigure}{0.8\textwidth}
		\centering
		\includegraphics[width=\textwidth]{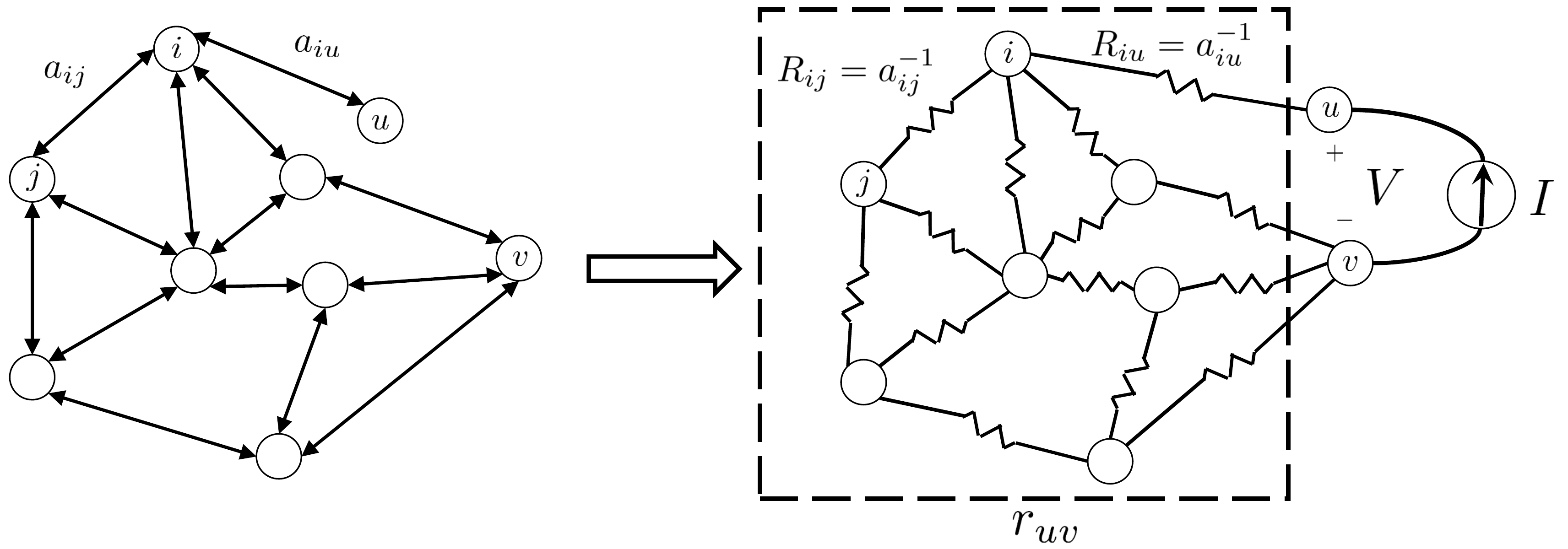}
		\caption{}
		\label{fig:fig3_1}
	\end{subfigure}%

	\centering
	\begin{subfigure}{0.8\textwidth}
		\centering
		\includegraphics[width=\textwidth]{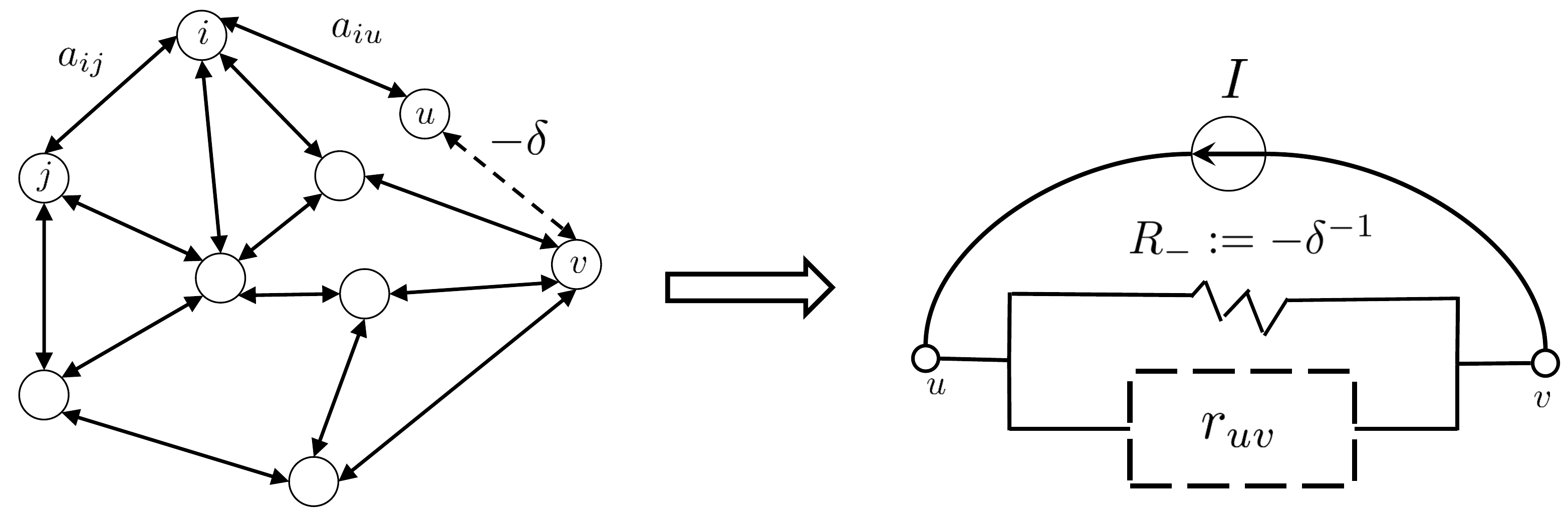}
		\caption{}
		\label{fig:fig3_2}
	\end{subfigure}%
	~ 
	\caption{(a) An equivalent electrical circuit for a connected undirected graph $\mathcal{G} (\mathcal{V},\mathcal{E}^+,\mathcal{W}^+)$ with non-negative weights to measure the effective resistance between two nodes $u$ and $v$. (b) An equivalent electrical circuit for $\mathcal{G} (\mathcal{V},\mathcal{E},\mathcal{W})$ which is obtained by adding a negative edge to $\mathcal{G} (\mathcal{V},\mathcal{E}^+,\mathcal{W}^+)$     .}
	\label{fig:fig3}
\end{figure}

We end this subsection by commenting on the relationship between the aforementioned results and the concept of effective resistance introduced in the literature e.g. see~\cite{zelazo2015robustness} and references there in. First, consider an undirected graph $\mathcal{G} (\mathcal{V},\mathcal{E}^+,\mathcal{W}^+)$ with non-negative weights and the Laplacian matrix $\bar{L}^+$. Assume that the graph is connected. We construct an electrical network from the graph $\mathcal{G} (\mathcal{V},\mathcal{E}^+,\mathcal{W}^+)$ by replacing each weighted edge $(i,j) $  with a resistor $R_{ij} = a_{ij}^{-1}$ as shown in Figure~\ref{fig:fig3_1}. Then, a constant current source $I$ is connected between nodes $u$ and $v$, and the voltage at its terminal is calculated. The effective resistance between nodes $u$ and $v$ is computed by $r_{uv} =  \frac{V}{I}$. Using electrical circuit theory, it was shown that $r_{uv}$ can be obtained by~\cite{klein1993resistance}
\begin{equation}
\label{eq:effective_resistance_undirected}
r_{uv} = \left( \mathbf{e}_u - \mathbf{e}_v \right)^T    Q^T \left(\bar{L}^+ \right)^{-1} Q \left( \mathbf{e}_u - \mathbf{e}_v \right),
\end{equation}
which has the same expression as~(\ref{eq:cond_single_weight_general_connected_postive_final}) with $q_{uv} = q_{vu} = 1$.
\newline
Now consider the case where a negative undirected edge is added between the nodes $u$ and $v$ to construct the undirected graph $\mathcal{G}(\mathcal{V},\mathcal{E}, \mathcal{W})$ with the Laplacian matrix $L$. The corresponding equivalent electrical circuit contains two parallel resistors $r_{uv}$ and $R_{-} = -\delta^{-1}$\footnote{Negative resistance has a practical meaning in electrical circuit theory as there exist electrical components with this property~(see~\url{https://en.wikipedia.org/wiki/Negative_resistance}).}, which correspond to the aforementioned effective resistance and the added negative edge, respectively~(Figure~\ref{fig:fig3_2}). In this case, the equivalent resistance between nodes $u$ and $v$ is calculated by~$$R_{th}=\frac{-\delta^{-1} r_{uv}}{r_{uv} - \delta^{-1}}.$$
The interpretation of the inequality~(\ref{eq:cond_single_weight_general_connected_postive_final})\footnote{Note that inequality~(\ref{eq:cond_single_weight_general_connected_postive_final}) is necessary and sufficient for this case as the underlying graph is undirected.} from the electrical circuit perspective is as follows: the equivalent resistance $R_{th}$ is positive as long as~(\ref{eq:cond_single_weight_general_connected_postive_final}) holds. Otherwise, $R_{th}$ is either short circuit~(if $\delta^{-1} = r_{uv}$) or negative~(if $\delta^{-1} < r_{uv}$). By taking into account this point and Corollary~\ref{corollary:one_negative_weight_undirected}, it is concluded that the Laplacian matrix $L$ satisfy the condition~(\ref{eq:desired_eigenvalues_condition}) if and only if the equivalent resistance $R_{th}$ of the corresponding equivalent electrical circuit is positive.    %

For directed graphs satisfying in the conditions of Corollary~\ref{corollary:one_negative_weight_undirected_1}, one may be interested in interpreting the condition~(\ref{eq:cond_single_weight_general_connected_postive_final}) using the concept of effective resistance similar to undirected graphs if directed graph satisfy assumptions of Corollary~\ref{corollary:one_negative_weight_undirected_1}. However, this is not directly possible. Indeed, the notion of effective resistance has been recently introduced for both directed and undirected graphs as~\cite{young2015new},
\begin{equation}
\label{eq:effective_resistance_directed}
r_{uv} = 2 \left( \mathbf{e}_u - \mathbf{e}_v \right)^T   Q^T \Sigma \, Q \left( \mathbf{e}_u - \mathbf{e}_v \right),
\end{equation}
where $\Sigma$ is a symmetric matrix obtained from the Lyapunov equation
\begin{equation}
\label{eq:Lyapunov_resistance}
\begin{aligned}
\bar{L} \Sigma &+ \Sigma \left( \bar{L} \right)^T = I_{N-1}.
\end{aligned}
\end{equation}
For a connected undirected graph $\mathcal{G} (\mathcal{V},\mathcal{E}^+,\mathcal{W}^+)$, $\frac{1}{2} \left(\bar{L}^+ \right)^{-1} = \frac{1}{2} \left(\bar{L}^{+} \right)^{-T}$ is a unique solution of~(\ref{eq:Lyapunov_resistance}) and $r_{uv}$ has the same expression as~(\ref{eq:effective_resistance_undirected}). Unlike undirected graphs, the expression of $r_{uv}$ in~(\ref{eq:cond_single_weight_general_connected_postive_final})~(even with $q_{uv} = q_{vu} = 1$) cannot be obtained from~(\ref{eq:effective_resistance_directed}) and
(\ref{eq:Lyapunov_resistance}), since  $\left(\bar{L}^+ \right)^{-1} \neq  \left(\bar{L}^{+} \right)^{-T}$.
This reveals that, unlike undirected graph, the upper bound of $\delta$ cannot be interpreted by effective resistance defined in~(\ref{eq:effective_resistance_directed}). As a result, it is not clear how to generalize the notion of effective resistance for directed graphs with non-negative weights by using electrical circuit theory.

\subsection{Directed graphs with multiple negative weights}
\label{subsec:multiple_negative_weights}
In this subsection, we examine the eigenvalues of the Laplacian matrix for graphs with multiple negative weights by stating the following question.

\begin{question}
	Given a graph~$\mathcal{G}_1(\mathcal{V},\mathcal{E}_1, \mathcal{W}_1)$ with non-negative edge weights, identify a set of pairs of nodes $\mathcal{E}_2 \subseteq \{ \mathcal{V} \times \mathcal{V} \} \backslash \mathcal{E}_1$, $\mathcal{E}_2 = \{ (u_1,v_1),\dots,(u_l,v_l) \}$, called sensitive pairs of nodes such that if a new graph is constructed from $\mathcal{G}_1$ by connecting those pairs of nodes with negative weights whose magnitudes are infinitesimal, the Laplacian matrix of the new graph has at least one eigenvalue with negative real part.
\end{question}

\begin{figure}[H]
	\centering
	\includegraphics[scale=0.5]{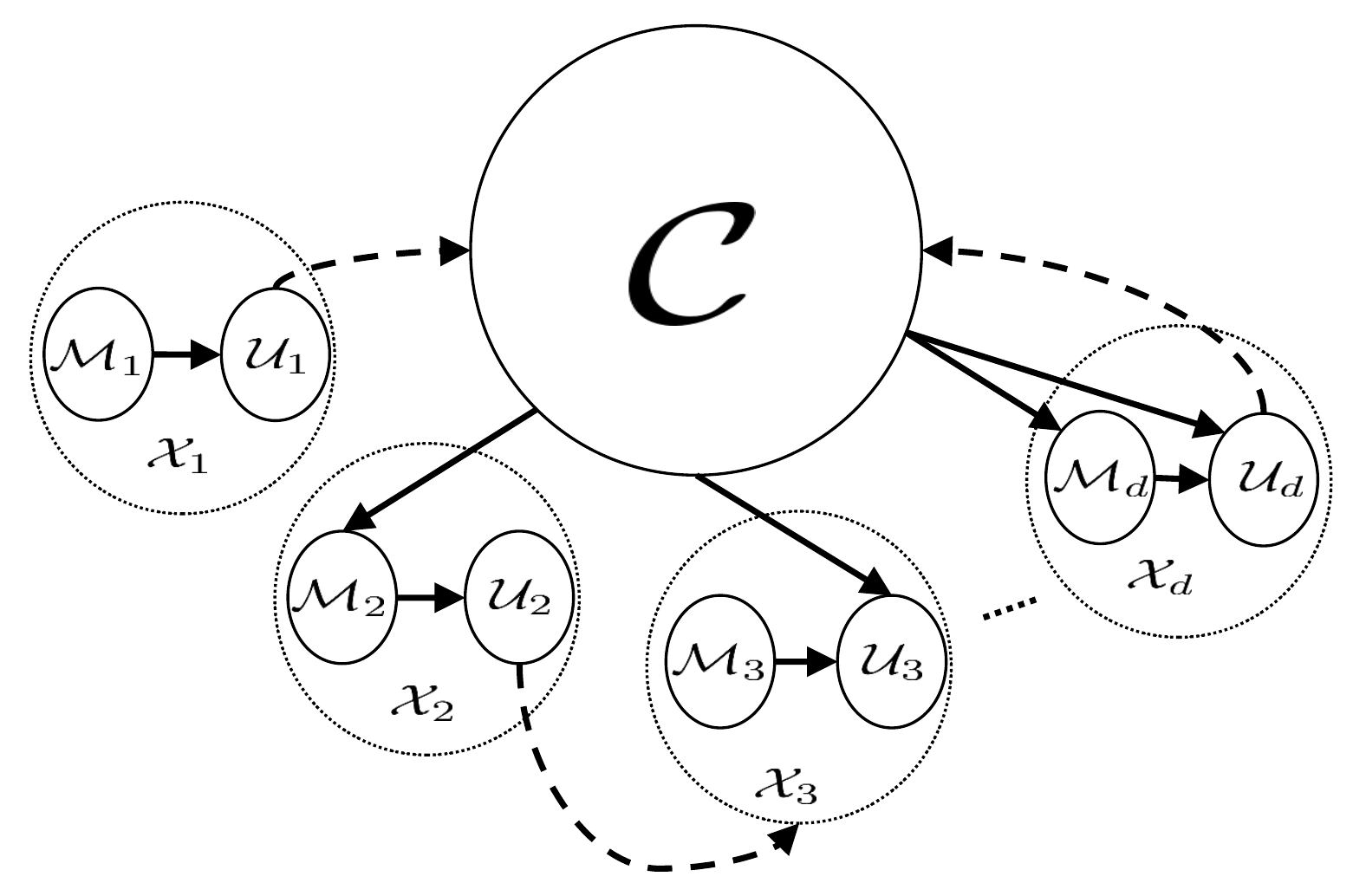}
	\caption{A schematic diagram of the connections between the nodes in a graph whose negative edges meet the conditions~$1$ and~$2$ in Theorem~\ref{theorem:connectivity_of_netowrk}. The solid arrows represent the positive weight edges while the dashed arrows show the edges with negative weights.}
	\label{fig:fig4}
\end{figure}

The following theorem is aimed at answering the aforementioned question. Figure~\ref{fig:fig4} shows a schematic diagram of a graph in which the negative edge weights connect those pairs of the nodes.

\begin{theorem}
	\label{theorem:connectivity_of_netowrk}
	Consider a graph~$\mathcal{G}_1(\mathcal{V},\mathcal{E}_1, \mathcal{W}_1)$ with non-negative edge weights and its Laplacian matrix $L_1$. Assume that $\mathcal{G}_1$ consists of $d \neq 1$ reach sets $\mathcal{R}_k,~k=1,\dots,d$ with corresponding reaching nodes sets~$\mathcal{U}_k$, exclusive sets~$\mathcal{X}_k$ and common sets~$\mathcal{C}_k$
	according to Definition~\ref{def:reachable_sets}.
	Let $\mathcal{G}(\mathcal{V},\mathcal{E}, \mathcal{W}) = \mathcal{G}_1(\mathcal{V},\mathcal{E}_1, \mathcal{W}_1) \bigoplus \mathcal{G}_2(\mathcal{V},\mathcal{E}_2, \mathcal{W}_2)$ where
	$\mathcal{E}_2 \subseteq \{ \mathcal{V} \times \mathcal{V} \} \backslash \mathcal{E}_1$, $\mathcal{E}_2 = \{ (u_1,v_1),\dots,(u_l,v_l) \}$, $\mathcal{W}_2(u_k,v_k) = - \epsilon a_{u_kv_k} $ for every $(u_k,v_k) \in \mathcal{E}_2$ with $\epsilon$ being a sufficiently small positive scalar and $a_{u_kv_k} >0$, and $\mathcal{W}_2(u_k,v_k) = 0 $ for every $(u_k,v_k) \in \{ \mathcal{V} \times \mathcal{V} \} \backslash \mathcal{E}_2$. If there exists at least one negative edge weight $(u_k,v_k) \in \mathcal{E}_2$ such that
	\begin{enumerate}
		\item{either $u_k \in \mathcal{U}_i$ and $v_k \in \mathcal{X}_j$ for some $j \neq i$,
		}
		\item {or $u_k \in \mathcal{U}_i$ and $v_k \in \mathcal{C}_j$ for some $j$ (including $j=i$),
		}
	\end{enumerate}
	then the Laplacian matrix of $\mathcal{G}$ has at least one eigenvalue with non-positive real part.
\end{theorem}

\textit{Proof.}
Suppose all conditions of Theorem~\ref{theorem:connectivity_of_netowrk} hold. The Laplacian matrix of $\mathcal{G}$ can be expressed as
\begin{equation}
	\label{eq:Laplacian_multiple_negative}
	\begin{aligned}
		L & = L_1 + L_2 = L_1 + \epsilon \tilde{L}_2,
	\end{aligned}
\end{equation}
where $L_2$ is the Laplacian matrix of $\mathcal{G}_2$, and $\tilde{L}_2$ is the Laplacian matrix for a graph $\widetilde{\mathcal{G}}=(\mathcal{V},\mathcal{E}_2,\widetilde{\mathcal{W}}_2)$ with $\widetilde{\mathcal{W}}_2(u_k,v_k) = - a_{u_kv_k} $ for every $(u_k,v_k) \in \mathcal{E}_2$, and $\widetilde{\mathcal{W}}_2(u_k,v_k) = 0 $ for every $(u_k,v_k) \in \{ \mathcal{V} \times \mathcal{V} \} \backslash \mathcal{E}_2$.

Without loss of generality, we suppose the nodes of $\mathcal{G}_1$ are labeled such that the structure of the adjacency matrix of $\mathcal{G}_1$ is consistent with the structure of~(\ref{eq:adjacency_general}).
Since zero is a semisimple eigenvalue of $L$ with multiplicity $d$, in view of Lemmas~\ref{lemma:perturbation} and \ref{lemma:Theta_diagonal_terms},
it is enough show that at least the real part of one of the eigenvalues of $\Theta := \Upsilon \tilde{L}_2 \Gamma$ becomes negative~(see Lemma~\ref{lemma:Theta_diagonal_terms} for definitions of~$ \Upsilon $ and $\Gamma$).
Conditions~$1$ and $2$ in the theorem statement are equivalent to statements~$3$ and~$4$ of Lemma~\ref{lemma:Theta_diagonal_terms}.
Hence, if there exists at least one negative edge $(u,v)$ satisfying  one of conditions~$1$ and~$2$, then according to Lemma~\ref{lemma:Theta_diagonal_terms}, $\Theta$ has at least a negative diagonal term  which means $\text{trace}(\Theta) < 0$ or equivalently $\sum\limits_{i=1}^{d} \lambda_i(\Theta) < 0$.
This implies that there exists at least one eigenvalue $\lambda_i(\Theta)$ with $\Re \{\lambda_i(\Theta) \} <0$. \qed

\begin{remark}
	For undirected graphs, it is straightforward to see that $\mathcal{R}_i = \mathcal{U}_i$. In view of Theorem~\ref{theorem:connectivity_of_netowrk}, if $\mathcal{G}$ is connected while $\mathcal{G}_1$ is not, then $\mathcal{G}_2$ contains negative weight edges which satisfy the first condition of Theorem~\ref{theorem:connectivity_of_netowrk}, and connect disconnected components of $\mathcal{G}_1$. In this case, Theorem~\ref{theorem:connectivity_of_netowrk} ensures that the Laplacian matrix of $\mathcal{G}$ is indefinite as it is symmetric and also has negative eigenvalue(s). This result is consistent with the result of Theorem~IV.3 in~\cite{zelazo2015robustness}.
\end{remark}

\begin{remark}
	\label{remark:more_cases}
	In Theorem~\ref{theorem:connectivity_of_netowrk}, we have characterized a class of sensitive pairs of nodes. However, there are other classes of sensitive pairs of nodes that we have not considered. As an example, $(u_k,v_k) \in \mathcal{E}_-$ such that:
	\begin{enumerate}
		\item {either $u \in \mathcal{U}_i~\text{and}~v \in \mathcal{X}_i,$
		}
		\item{ or $u \in \mathcal{R}_i \backslash \mathcal{U}_i~\text{and}~ v \in \mathcal{R}_j~\text{for some}~i \neq j.$
		}
	\end{enumerate}
	In this case, it can be shown that the corresponding diagonal term of the matrix $\Theta$ in the proof of Theorem~\ref{theorem:connectivity_of_netowrk} is zero. Hence, the same argument in Theorem~\ref{theorem:connectivity_of_netowrk} cannot be followed to achieve similar results.
\end{remark}

\begin{figure}[ht]
	\centering
	\includegraphics[scale=0.4]{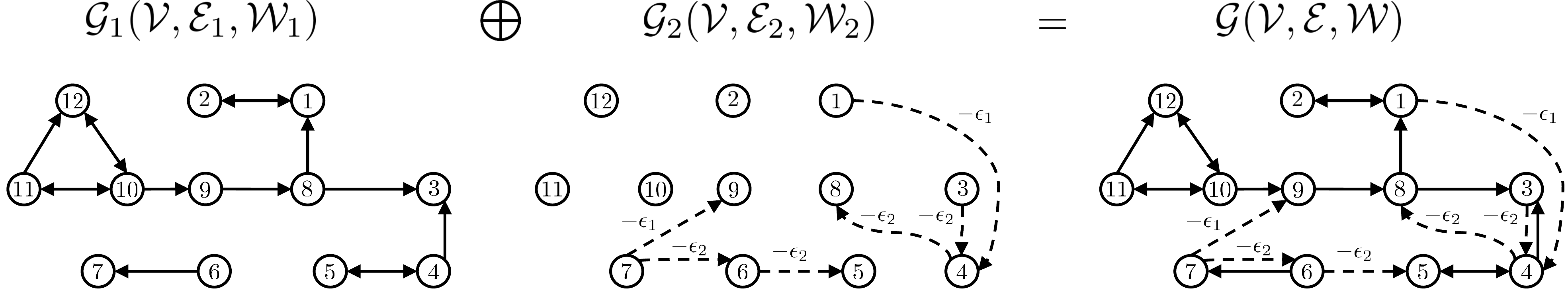}
	\caption{A directed signed graph studied in Example 3. The solid arrows represent the positive weight edges while the dashed arrows show the edges with negative weights. For the first case and second case $\epsilon_1 = 0.001, \epsilon_2 = 0$ and $\epsilon_1 = 0, \epsilon_2 = 0.001$, respectively.}
	\label{fig:fig5}
\end{figure}

\begin{example}
	We illustrate the results of Theorem~\ref{theorem:connectivity_of_netowrk} and the point mentioned in Remark~\ref{remark:more_cases} using the graphs depicted in Figure~\ref{fig:fig5}. All positive weights are assigned to be $2$.  We have observed in Example~1 that the graph $\mathcal{G}_1$ has three reaching nodes sets $\mathcal{U}_1 = \{ 1,2 \}$, $\mathcal{U}_2 = \{ 3 \}$, $\mathcal{U}_3 = \{ 7 \}$, three exclusive sets $\mathcal{X}_1 = \{ 1,2\}$, $\mathcal{X}_2 = \{ 3 ,4 ,5 \}$, $\mathcal{X}_3 = \{ 6, 7 \}$ , and three common sets $\mathcal{C}_1 = \mathcal{C}_2 = \{ 8,9,10,11,12 \}$ and $\mathcal{C}_3 = \emptyset$.

	In what follow, we demonstrate the application of Theorem~\ref{theorem:connectivity_of_netowrk}. To this end, let $\epsilon_1 = 0.0001$ and $\epsilon_2 = 0$. In this case, $\mathcal{E}_2 = \{(7,9),(1,4)\}$. Note that
	\begin{enumerate}
		\item[(i)] the nodes $1$ and $4$ belong to the sets $\mathcal{U}_1$ and $\mathcal{X}_2$~(Condition 1 in Theorem~\ref{theorem:connectivity_of_netowrk});
		\item[(ii)] the nodes $7$ and $9$ belong to the sets $\mathcal{U}_3$ and $\mathcal{C}_1( \text{or}~\mathcal{C}_2)$~(Condition 2 in Theorem~\ref{theorem:connectivity_of_netowrk}),
	\end{enumerate}
	which means that the negative edge satisfies the conditions of Theorem~\ref{theorem:connectivity_of_netowrk}. By computing the eigenvalues of the Laplacian matrix of $\mathcal{G}$, it is observed that it has two eigenvalues with negative real parts, while the remaining eigenvalues have positive real parts.

	Next, we investigate the statement of Remark~\ref{remark:more_cases} via choosing $\epsilon_1 = 0$, $\epsilon_2 = 0.0001$ and $\mathcal{E}_2 =  \{(3,4), (7,6),(6,5),(4,8)\}$. Note that
	\begin{enumerate}
		\item[(i)] the nodes $3$ and $4$ belong to the sets $\mathcal{U}_2$ and $\mathcal{X}_2$~(Condition 1 in Remark~\ref{remark:more_cases});
		\item[(ii)] the nodes $7$ and $6$ belong to the sets $\mathcal{U}_3$ and $\mathcal{X}_3$~(Condition 1 in Remark~\ref{remark:more_cases});
		\item[(iii)] the nodes $6$ and $5$ belong to the sets $\mathcal{X}_3 \subset \mathcal{R}_3 \backslash \mathcal{U}_3$ and $\mathcal{X}_2 \subset \mathcal{R}_2$~(Condition 2 in Remark~\ref{remark:more_cases});
		\item[(iv)] the nodes $4$ and $8$ belong to the sets $\mathcal{X}_2 \subset \mathcal{R}_2 \backslash \mathcal{U}_2$ and $\mathcal{C}_1 \subset \mathcal{R}_1$~(Condition 2 in Remark~\ref{remark:more_cases}),
	\end{enumerate}
	which means that the negative edges of the graph depicted in Figure~\ref{fig:fig5} meet the conditions outlined in Remark~\ref{remark:more_cases}.
	By computing the eigenvalues of the Laplacian matrix, we noticed that all eigenvalues of the Laplacian matrix (excluding the single zero one) have positive real parts, however, their real parts become negative if the absolute values of negative weights are not small. This shows that for these negative edges, knowing that they are negative is not enough to conclude about the sign of real parts of eigenvalues, but their weights also play a role, unlike the one in the first case.
\end{example}

\section{Application}
\label{sec:application}
In this section, we elucidate the application of our results to the consensus problem in social networks with antagonistic interactions. The dynamic of the whole network with $N$ node is written as
\begin{equation}
\dot{\mathbf{x}}= - L \mathbf{x},
\end{equation}
where $\mathbf{x} \in \mathbb{R}^N,
~\mathbf{x}=[x_1 x_2 \dots x_N]^T$ is the stack vector consisting of the state of each node in the network. $L$ is the Laplacian matrix of the network graph. The network achieves consensus, i.e. $\lim\limits_{t \to \infty} \| x_i(t) - x_j(t) \| = 0$ for all $i,j =1,\dots,N$ and $ \mathbf{x}(0) \in \mathbb{R}^N$, if and only if the condition~(\ref{eq:desired_eigenvalues_condition}) holds ~\cite[Lemma~2]{olfati2007consensus}.

Consider a network with a graph network $\mathcal{G}_1(\mathcal{V},\mathcal{E}_1,\mathcal{W}_1)$ in Example~$1$. In Example~$1$, we have shown that the Laplacian matrix $L_1$  with these values for weights satisfies~(\ref{eq:desired_eigenvalues_condition}), and consequently, the network achieves consensus. We have perturbed the network with negative directed edge $a_{38} = -\delta$. We have pointed  out that  the eigenvalues of $L$  meet~(\ref{eq:desired_eigenvalues_condition}) if and only if $\delta < 1.94285 $.

\begin{figure}[h!]

	\centering
	\begin{subfigure}{0.5\textwidth}
		\centering
		\includegraphics[width=\textwidth]{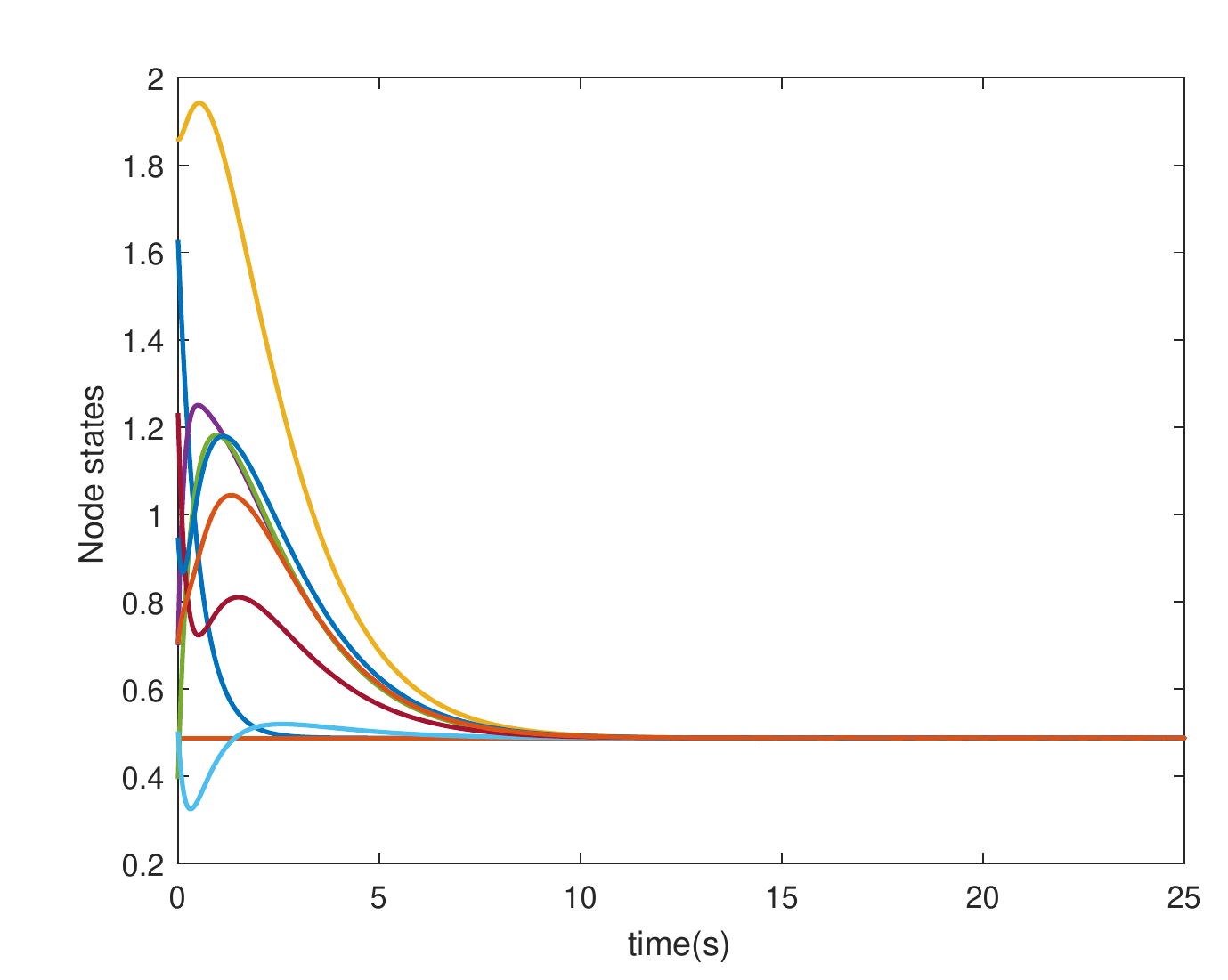}
		\caption{}
		\label{fig:fig6_2}
	\end{subfigure}%
	\begin{subfigure}{0.5\textwidth}
		\centering
		\includegraphics[width=\textwidth]{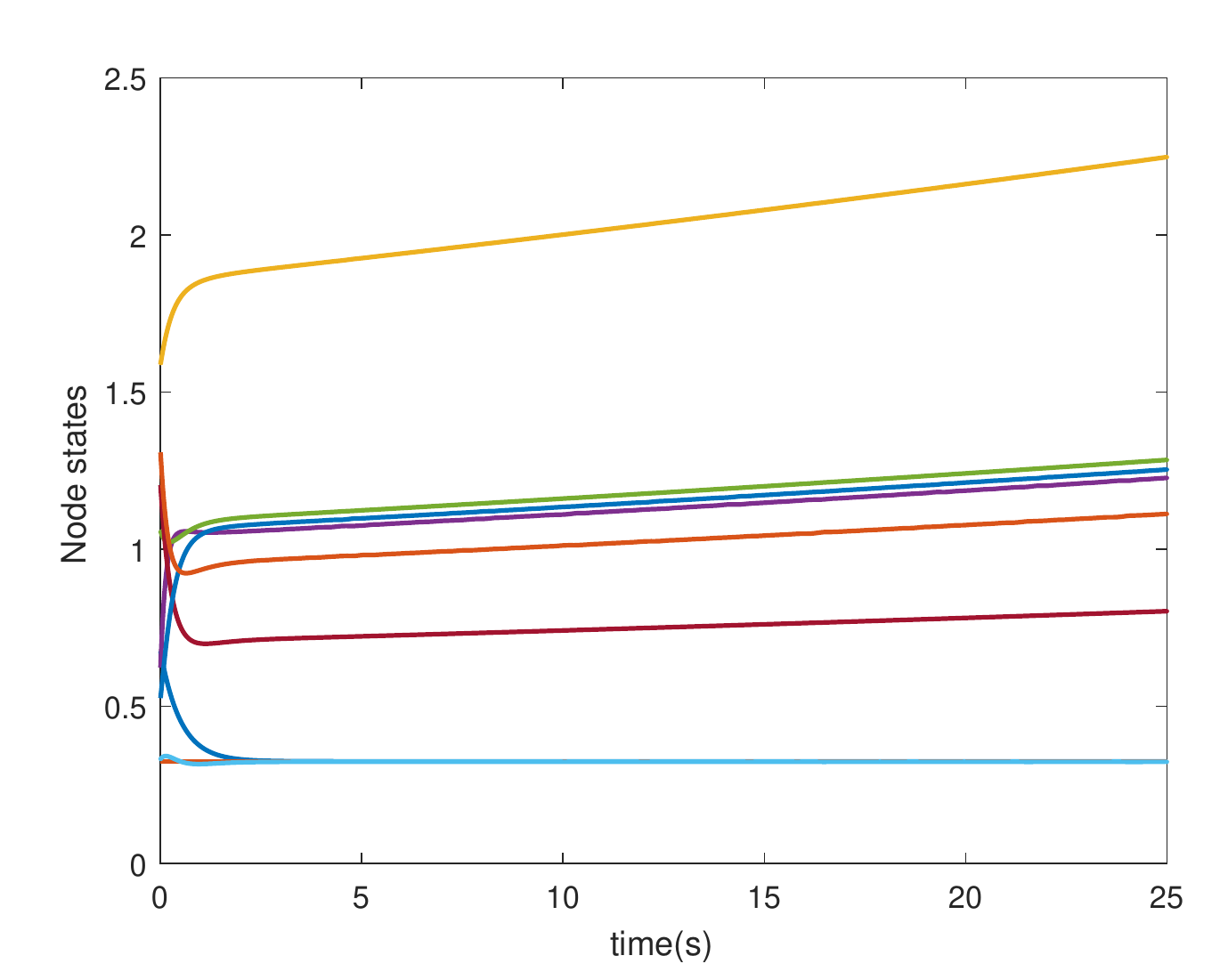}
		\caption{}
		\label{fig:fig6_3}
	\end{subfigure}%
	~ 
	\caption{ Evolution of node states over time for a network with the graph in Example~$1$. The perturbation between two nodes $3$ and $8$ set equal to (a) $\delta_{38} = 1.5$ (b)~ $\delta_{38} = 1.95$.}
	\label{fig:fig6}
\end{figure}

With the choice of $ \delta= 1.5$, the eigenvalues of $L$ satisfy~(\ref{eq:desired_eigenvalues_condition}) and network achieves consensus as shown in Figure~\ref{fig:fig6_2}. However, with $\delta= 1.95$, the eigenvalues of $L$ no longer satisfy~(\ref{eq:desired_eigenvalues_condition}) and network does not achieve consensus as shown in Figure~\ref{fig:fig6_3}.

\section{Conclusion}
We have analyzed the eigenvalues of the Laplacian matrix for a class of directed graphs with positive and negative weights. The first part of our paper has dealt with directed signed graphs where  arbitrary negative edges were added between any arbitrary pairs of nodes in the original graphs. Under certain conditions, we have shown that if the magnitude of that added negative weight violates a computable upper bound, all eigenvalues of the Laplacian matrix (except the single zero eigenvalue) of the original graph have positive real parts. We have pointed out that, unlike undirected graphs, it appears that effective resistance definition does not apply directly to this case. In the second part, we have shown that if the graph is perturbed by adding edge(s) with negative weight(s) to certain nodes, the Laplacian matrix of the perturbed graph has at least one eigenvalue with negative real part.

\section*{Acknowledgement}

The error in the original paper was pointed out by Dr. Sei Zhen Khnog from University of Minnesota, as well as Dr. Daniel Zelazo and Dr. Dwaipayan Mukherjee from Israel Institute of Technology. We appreciate their insightful skepticism which enabled us to correct our results.
\appendix

%

%
\section{}
\label{Appendix:a}
The following lemmas are used in the proofs of Theorems~\ref{theorem:one_negative_weight_undirected_necessity} and~\ref{theorem:one_negative_weight_undirected_corrected}.

\begin{lemma}[Fact 2.16.3 in~\cite{bernstein2009matrix}]
	\label{lemma:rank_one_updates}
	Let $A \in \mathbb{C}^{n \times n}$, assume that $A$ is nonsingular, and let $c,d \in \mathbb{R}^{n \times 1}$. Then
	\begin{equation}
	\label{eq:rank_one_updates}
	\det (A+cd^T) = \det(A) (1+d^T A^{-1} c).
	\end{equation}
\end{lemma}

\begin{lemma}
	\label{lemma:one_negative_weight_undirected}
	Consider a signed graph~$\mathcal{G}_1(\mathcal{V},\mathcal{E}_1, \mathcal{W}_1)$ with the Laplacian matrix $L_1$. Assume that $L_1$ has only one zero eigenvalue and the rest of its eigenvalues have positive real parts.
	Construct a new graph $\mathcal{G}(\mathcal{V},\mathcal{E}, \mathcal{W}) = \mathcal{G}_1(\mathcal{V},\mathcal{E}_1, \mathcal{W}_1) \bigoplus \mathcal{G}_2(\mathcal{V},\mathcal{E}_2, \mathcal{W}_2)$,
	where $ \mathcal{E}_2= \{ (u,v), (v,u)\}$ and $\mathcal{W}_2(u,v) = -\delta_{uv} ,~\mathcal{W}_2(v,u)= -\delta_{vu}$ with $\delta_{uv} \geq 0,~\delta_{vu} \geq 0 $. Denote $L$ the Laplacian matrix of $\mathcal{G}$. Then,
	\begin{equation}
		\label{eq:spectrum_L_1bar_L_bar}
		\text{Spec} \{ \bar{L}_1^{-1} \bar{L} \} = \{ \underbrace{ 1,\dots, 1}_{N-2}, 1-r(\omega,\delta_{uv},\delta_{vu}) \},
	\end{equation}
	where $r(\omega,\delta_{uv},\delta_{vu}) = (\mathbf{e}_u - \mathbf{e}_v)^T Q^T \bar{L}_1^{-1} Q (\delta_{uv} \mathbf{e}_u - \delta_{vu} \mathbf{e}_v)$, and $\bar{L} = Q L Q^T$ and $\bar{L}_1 = Q L_1 Q^T$ the reduced Laplacian matrices for $\mathcal{G}$ and $\mathcal{G}_1$, respectively.
\end{lemma}

\textit{Proof.}
	Denote $L_2$ the Laplacian matrix of $\mathcal{G}_2$. Since $\mathcal{E}_2= \{ (u,v), (v,u)\}$, $L_2$ can be expressed as $L_2 = - (\delta_{uv} \mathbf{e}_u - \delta_{vu} \mathbf{e}_v) (\mathbf{e}_u - \mathbf{e}_v)^T$. Furthermore,
	\begin{equation}
		\label{eq:spec_L_2}
		\text{Spec}(L_2) = \{ \underbrace{ 0,\dots, 0}_{N-1}, -(\delta_{uv} + \delta_{vu})   \},
	\end{equation}
	with the set of eigenvectors $ \{ \mathbf{1}_N, \mathbf{e}_1,\dots,\mathbf{e}_{N} \} \backslash \{{\mathbf{e}_u,\mathbf{e}_v} \} $ that corresponds to the zero eigenvalues, and $\delta_{uv} \mathbf{e}_u - \delta_{vu} \mathbf{e}_v$ that corresponds to the eigenvalue $-(\delta_{uv} + \delta_{vu})$. Using the first property in Lemma~\ref{lemma:eigenvalues_relation} with~(\ref{eq:spec_L_2}), we have
	\begin{equation}
		\label{eq:spec_L_bar_2}
		\text{Spec} \left( \bar{L}_2 \right) = \{ \underbrace{ 0,\dots, 0}_{N-2}, -(\delta_{uv} + \delta_{vu})   \},
	\end{equation}
	where $\bar{L}_2 = Q L_2 Q^T$. The graphs $\mathcal{G}_1$ and $\mathcal{G}_2$ have the same set of nodes which means  $L$ can be written as $L=L_1 + L_2$  leading to the following expression for $\bar{L}$,
	\begin{equation}
		\label{eq:Lbar_new_representation}
		\bar{L} =\underbrace{Q L_1 Q^T}_{\bar{L}_1} - Q (\delta_{uv} \mathbf{e}_u - \delta_{vu} \mathbf{e}_v)  (\mathbf{e}_u - \mathbf{e}_v)^T Q^T.
	\end{equation}
	Since $L_1$ has only one zero eigenvalue, $\bar{L}_1$ is invertible according to Lemma~\ref{lemma:eigenvalues_relation}. By multiplying both sides of~(\ref{eq:Lbar_new_representation}) by $\bar{L}_1^{-1} $ and then $(\mathbf{e}_u - \mathbf{e}_v)^T Q^T$, we obtain
	\begin{equation}
		\label{eq:eigenvalue_necessity}
		\begin{aligned}
			(\mathbf{e}_u - \mathbf{e}_v)^T Q^T \bar{L}_1^{-1} \bar{L} &= (\mathbf{e}_u - \mathbf{e}_v)^T Q^T- \underbrace{(\mathbf{e}_u - \mathbf{e}_v)^T Q^T \bar{L}_1^{-1} Q (\delta_{uv} \mathbf{e}_u - \delta_{vu} \mathbf{e}_v)}_{r_{\delta} := r(\omega,\delta_{uv},\delta_{vu})}  (\mathbf{e}_u - \mathbf{e}_v)^T Q^T\\
			&= (\mathbf{e}_u - \mathbf{e}_v)^T Q^T (1- r_{\delta}),
		\end{aligned}
	\end{equation}
	or equivalently,
	\begin{equation}
		\label{eq:eigenvalue_necessity_new_presentation}
		(\mathbf{e}_u - \mathbf{e}_v)^T Q^T \left( \bar{L}_1^{-1} \bar{L} - (1- r_{\delta}) I_{N-1} \right)= \mathbf{0}_{N-1}^T.
	\end{equation}
	Since $Q^T$ is a full column rank matrix, (\ref{eq:eigenvalue_necessity_new_presentation}) implies that the vector $Q E_2\in \mathbb{R}^{ N-1}$ is a left eigenvector of matrix $\bar{L}_1^{-1} \bar{L} \in \mathbb{R}^{(N-1) \times (N-1)}$ that corresponds to the eigenvalue $1- r_{\delta}$. Showing~(\ref{eq:spectrum_L_1bar_L_bar}) is equivalent to showing that $X= I_{N-1} - \bar{L}_1^{-1} \bar{L} $ has $N-2$ zero eigenvalues.  From~(\ref{eq:Lbar_new_representation}), we obtain
	\begin{equation}
		\begin{aligned}
			\bar{L}_1 X & = \bar{L}_1 - \bar{L} = - \bar{L}_2,
		\end{aligned}
	\end{equation}
	which means $\text{spec} \{ \bar{L}_1 X \} = \text{spec} \{ - \bar{L}_2 \}$. Using~(\ref{eq:spec_L_bar_2}) and noting that $\bar{L}_1$ is non-singular, we conclude $X$ has $N-2$ zero eigenvalues. This completes the proof. \qed
	%

\section{}
\label{Appendix:b}
We recall particular case of Theorem~$2.1$ in~\cite{Moro1997lidskii} and present Lemma~\ref{lemma:Theta_diagonal_terms} that are used in proof of~Theorem~\ref{theorem:connectivity_of_netowrk}.

\begin{lemma}[\cite{Moro1997lidskii}]
	\label{lemma:perturbation}
	Assume $\lambda \in \mathbb{C}$ is semisimple eigenvalue of a square matrix $A \in \mathbb{R}^{N \times N}$ with multiplicity $d$. Consider a perturbed matrix $A + \epsilon B$ where $B$ is an arbitrary matrix and $\epsilon$ introduces a small perturbation parameter. Then, there are $d$ eigenvalues of the perturbed matrix which are described by a first-order expansion
	\begin{equation}
		\lambda_i=\lambda + \xi_i \epsilon + o(\epsilon),~i=1,\dots,d
	\end{equation}
	where $\xi_i$ are the eigenvalues of the $d \times d$ matrix $ \Upsilon B \Gamma$ and $o(\epsilon)$ contains remaining terms such that $\lim\limits_{\epsilon \to 0} \frac{o(\epsilon)}{\epsilon} = 0$. The $j^{th}$ row and column of $\Upsilon$ and $\Gamma$ are respectively, the left and right eigenvectors of $A$ corresponding to $\lambda$ which are orthonormal, i.e. $\Upsilon \Gamma = I_d$.
\end{lemma}

\begin{lemma}
	\label{lemma:Theta_diagonal_terms}
	Consider a graph~$\mathcal{G}_1(\mathcal{V},\mathcal{E}_1, \mathcal{W}_1)$ with non-negative edge weights and its Laplacian matrix $L_1$. Assume,
	\begin{enumerate}
		\item {$\mathcal{G}_1$ consists of $d \neq 1$ reach sets $\mathcal{R}_k,~k=1,\dots,d$.
		}
		\item{All nodes of $\mathcal{G}_1$ are labeled such that the structure of the adjacency matrix of $\mathcal{G}_1$ has the structure in~(\ref{eq:adjacency_general}).
		}
	\end{enumerate}
	Denote $\gamma_k$ and $\mu_k $  the right and left eigenvectors associated with the zero eigenvalue of $L_1$, respectively. Construct a new graph $\mathcal{G}(\mathcal{V},\mathcal{E}, \mathcal{W}) = \mathcal{G}_1(\mathcal{V},\mathcal{E}_1, \mathcal{W}_1) \bigoplus \mathcal{G}_2(\mathcal{V},\mathcal{E}_2, \mathcal{W}_2)$ with $ \mathcal{E}_2 \subseteq \{ \mathcal{V} \times \mathcal{V} \} \backslash \mathcal{E}_1 $, and $\mathcal{W}_2(u,v) < 0 $ for every $(u,v) \in \mathcal{E}_2$. Define $\Theta \in \mathbb{R}^{d \times d}$ as
	\begin{equation}
		\label{eq:definition_Theta}
		\Theta  = \Upsilon L_2 \Gamma,
	\end{equation}
	where $\Upsilon = \left[ \frac{\mu_1}{\| \mu_1 \|} \dots \frac{\mu_d}{ \| \mu_d \|} \right]^T$, $\Gamma = \left[ \frac{\gamma_1}{ \| \gamma_1 \|} \dots \frac{\gamma_d}{ \| \gamma_d \|} \right]$, and $L_2$ is the Laplacian matrix of $\mathcal{G}_2$. If $| \mathcal{E}_1 | = 1$, the following statements are true for $i=1,\dots,d$.
	\begin{enumerate}
		\item{ If $u \notin \mathcal{U}_i$, then $[\Theta]_{ii}= 0$;
		}
		\item{ if $u \in \mathcal{U}_i$ and $v \in \mathcal{X}_i$, then $[\Theta]_{ii}= 0$;
		}
		\item{ if $u \in \mathcal{U}_i$ and $ v \in \mathcal{X}_j$ or  $ v \in \mathcal{C}_j$ for $ j \neq i$, then  $[\Theta]_{ii} < 0$;
		}
		\item{ if $u \in \mathcal{U}_i$ and $ v \in \mathcal{C}_i$, then $[\Theta]_{ii} < 0$,
		}
	\end{enumerate}
	where the reaching nodes sets~$\mathcal{U}_k$, the exclusive sets~$\mathcal{X}_k$ and the common sets~$\mathcal{C}_k$ are described according to Definition~\ref{def:reachable_sets}.

	Furthermore, if $| \mathcal{E}_2 | > 1$, then $[\Theta]_{ii} \leq 0$. In this case,  $[\Theta]_{ii} < 0$ if and only if there exists at least one edge $(u,v) \in \mathcal{E}_2$ which satisfies either of the conditions given in statements $3$ or $4$.
\end{lemma}

\textit{Proof.}
Suppose all conditions of Lemma~\ref{lemma:Theta_diagonal_terms} hold. Partition the matrix $L_2$ in the same way as the adjacency matrix of $\mathcal{G}_1$. Since all weights of $\mathcal{G}_1$ are non-negative,  $\gamma_k$ and $\mu_k$ are characterized according to Lemma~\ref{lemma:eigenvectors_properties}. Taking into account the structure of $\mu_i$, we achieve
$$
[\Theta]_{ii} =  \tilde{\mu}^T_i \underbrace{\left[\begin{array}{cccccc}
	L_{2_{\mathcal{U}_i \mathcal{X}_1}}&L_{2_{\mathcal{U}_i \mathcal{X}_2}}&\hdots&L_{2_{\mathcal{U}_i \mathcal{X}_d}}&L_{2_{\mathcal{U}_i \mathcal{C}}}
	\end{array}\right]}_{L_{2_{\mathcal{U}_i}}} \frac{\gamma_i}{\|\gamma_i \|},
$$
where $ L_{2_{\mathcal{U}_i \mathcal{X}_j}}= \left[\begin{array}{cc}
L_{2_{\mathcal{U}_i \mathcal{U}_j}}&L_{2_{\mathcal{U}_i \mathcal{M}_j}}
\end{array}\right]$, and $\tilde{\mu}_i \in \mathbb{R}^{|\mathcal{U}_i|}$ contains the nonzero terms of the vector $\frac{\mu_i}{ \| \mu_i \|}$. First, assume there is only a single negative edge $\mathcal{E}_2= \{ (u,v) \}$ with $\mathcal{W}_2(u,v) = - \delta < 0$.

Proof of Statement~1: If $u \notin \mathcal{U}_i$, then $L_{2_{\mathcal{U}_i}} = \mathbf{0}$ leading to $[\Theta]_{ii} = 0$.

Proof of Statement~2: If $u \in \mathcal{U}_i$, the expression $[\Theta]_{ii}$ in above can be expanded as,
\begin{equation}
	\label{eq:big_theta_diagonal_term}
	\begin{aligned}
		~[\Theta]_{ii} &=-\frac{\delta}{ \| \mu_i \| \| \gamma_i \|} [\mu_i]_u [\gamma_i]_u + \frac{\delta}{ \| \mu_i \|  \| \gamma_i \|} [\mu_i]_u [\gamma_i]_v \\
		&= -\frac{\delta}{ \| \mu_i \| \| \gamma_i \|} [\mu_i]_u \left( 1 - [\gamma_i]_v \right).
	\end{aligned}
\end{equation}
According to Lemma~\ref{lemma:eigenvectors_properties}, $[\mu_i]_u > 0$ and $[\gamma_i]_v =1$ for $v \in \mathcal{X}_i$ which along with~(\ref{eq:big_theta_diagonal_term}) concludes $[\Theta_{ii}] = 0$.

Proof of Statements 3 and 4: For $v$ given in statement 3 or statement 4,  $[\gamma_i]_v =0$ or $[\gamma_i]_v <1$ which using~(\ref{eq:big_theta_diagonal_term}) show $[\Theta]_{ii} < 0$.

To deal with multiple negative edge weights, i.e. $| \mathcal{E}_2 | > 1$, it should be noted that the Laplacian matrix $L_2$ can be written as a sum of Laplacian matrices each of which corresponding to one negative edge weight. As adding a negative edge weight does not contribute to any positive value for  $[\Theta]_{ii}$, we conclude that the diagonal terms of $\Theta$ are non-positive, i.e. $[\Theta]_{ii} \leq 0$ for $i=1,\dots,d$. In addition,  from our discussion for a single negative edge weight, it is observed that $[\Theta]_{ii} $ becomes negative if and only if there exists at least one edge $(u,v) \in \mathcal{E}_2$ that meets the one of the conditions mentioned in the statements~$3$ and~$4$. This completes the proof. \qed
%



\bibliographystyle{elsarticle-harv}
\bibliography{bib_negative_directed}

\end{document}